%% file: slopeobs_archive.tex
\documentclass[11pt]{amsart} 
\usepackage[dvips]{graphicx,color}
\usepackage{amssymb,amscd,latexsym,epsfig}
\usepackage[arrow,matrix,line]{xy} \usepackage[mathscr]{euscript}
\DeclareFontFamily{OT1}{rsfs}{}
\DeclareFontShape{OT1}{rsfs}{n}{it}{<-> rsfs10}{}
\DeclareMathAlphabet{\curly}{OT1}{rsfs}{n}{it}

\newcommand{\rt}[1]{\stackrel{#1\,}{\rightarrow}}

\newcommand\comp{{\,}_{{}^\circ}} 
\newcommand\C{\mathbb C} \newcommand\Q{\mathbb Q}
\newcommand\R{\mathbb R} 
 
\newcommand\PP{\mathbb P} \renewcommand\P{\mathcal P}
 \newcommand\X{\curly X}
 
\renewcommand\L{\mathcal L} 
\renewcommand\O{\mathcal O} \newcommand\I{\curly I}
\newcommand\into{\hookrightarrow}

\renewcommand\_{^{\ }_} 
\newcommand{\rk}{\operatorname{rank}}

\newcommand{\Hom}{\operatorname{Hom}}

\newcommand{\vol}{\operatorname{vol}}

\newcommand{\ch}{\operatorname{ch}} \newcommand{\aut}{\mathfrak{aut}}

\makeatletter \@addtoreset{equation}{section} \makeatother

\newtheorem{thm}[equation]{Theorem} \newtheorem{lem}[equation]{Lemma}
\newtheorem{cor}[equation]{Corollary}
\newtheorem{prop}[equation]{Proposition}
\newtheorem{conj}[equation]{Conjecture}

\theoremstyle{definition} \newtheorem{defn}[equation]{Definition}
\newtheorem{rmk}[equation]{Remark}
\newtheorem{rmks}[equation]{Remarks}
\newtheorem{example}[equation]{Example}

\title[An obstruction to the existence of cscK metrics] {\textbf  An obstruction
to the existence
of constant scalar curvature K\"{a}hler metrics.}

\author[J. Ross and R. P. Thomas]{Julius Ross and Richard Thomas}

\begin{document}

\begin{abstract} \noindent
  We prove that polarised manifolds that admit a constant scalar curvature
  K\"ahler (cscK) metric satisfy a
  condition we call slope semistability. That is, we define the slope
  $\mu$ for a projective manifold and for each of its subschemes, and
  show that if $X$ is cscK then $\mu(Z)\le\mu(X)$ for all subschemes $Z$.
  
  This gives many examples of manifolds with K\"ahler classes which do not admit cscK metrics,
  such as del Pezzo surfaces and projective bundles.  If $\PP(E)\to B$ is
  a projective bundle
  which admits a cscK metric in a rational K\"ahler class with sufficiently
  small fibres, then $E$ is a slope semistable bundle (and $B$ is a slope
  semistable polarised manifold). The same is true for \emph{all} rational
  K\"ahler classes if the base $B$ is a curve.
  
  We also show that the slope inequality holds automatically for smooth
  curves, canonically polarised and Calabi Yau manifolds, and
  manifolds with $c_1(X)<0$ and $L$ close to the canonical
  polarisation.
\end{abstract}

\maketitle


\section{Introduction}

An important problem in K\"{a}hler geometry is that of finding a
constant scalar curvature K\"{a}hler (cscK) metric in a given K\"ahler
class on a complex manifold $X$. For a curve this is provided by the
uniformisation theorem. For general $X$ the class $[\omega]\in
H^2(X,\R)$ admits a K\"ahler-Einstein metric (which is therefore cscK)
when $c_1(X)=0$ \cite{Y1}, or when $c_1(X)<0$ and
$[\omega]=-\lambda[c_1(X)]$ \cite{Au, Y1}.

The first known obstructions to the existence of cscK metrics came
from the holomorphic automorphism group.  The most famous is the
Calabi-Futaki invariant of the K\"ahler class. This is a character on
the Lie algebra $\aut(X)$ of the holomorphic automorphism group which
must vanish if the class admits a cscK metric \cite{Fut}.

Tian defined a finer obstruction called K-stability, arising from
certain degenerations (or \emph{test configurations}) of $X$
\cite{Ti2, Ti3}. Moreover it is conjectured that K-polystability is a
necessary and sufficient condition for the existence of cscK metrics;
see Conjecture \ref{conj:cscK}. One direction of this conjecture is now almost
proved: it is known that cscK implies K-semistability \cite{Do5}. Thus test configurations can provide obstructions to cscK
metrics. In particular those arising from a $\C^\times$-action recover
the Calabi-Futaki obstruction, and these ``product configurations" are
currently the only test configurations that have been systematically
studied.

In this paper we consider test configurations canonically associated
to subschemes of $X$, yielding a new obstruction to the existence of
cscK metrics.  These configurations are more complicated than product
configurations, in particular the central fibre is non-normal. The
motivation is an analogy with stability for vector bundles; just as
subsheaves can destabilise a sheaf or bundle we show how subschemes
can destabilise $X$. In Section \ref{sec:definitionofslope} we define,
by analogy with vector bundles, a notion of slope (semi)stability of a
manifold and rational K\"ahler class.  We prove in Section
\ref{sec:kstabilityimpliesslopestability} that this gives an
obstruction to K-semistability, and hence to cscK metrics.

That is K-semistability implies slope semistability.  A partial
converse is given in Theorem 6.1 of \cite{RT}; in particular the two
are shown to be equivalent for curves. In trying to form moduli of
varieties in algebraic geometry using Geometric Invariant Theory,
other stability conditions arise, for instance Chow stability.
Different notions of slope are given for these in \cite{RT}, and
stability is shown to imply the relevant slope stability (Proposition
4.33 and Theorem 7.2 of \cite{RT}).

In Section \ref{sec:nonpositive} we show slope stability for the
canonical class when $c_1(X)<0$, and for arbitrary classes when
$c_1(X)=0$ (as expected from the existence of their K\"ahler-Einstein
metrics).  We also show slope stability for classes close to the
canonical class when $c_1(X)<0$, and compare to some similar
analytical results of \cite{We}.  Slope stability of smooth curves is
proved in \ref{sec:slopestabilitysmoothcurves}, which by Corollary 6.7
of \cite{RT} implies K-stability.  As far as we know this is the only
direct, non-analytic proof of K-stability of smooth curves.

We apply the slope formula to study unstable projective bundles in
Section \ref{proj}, providing a converse to the results of Hong \cite{Ho}.
When the base is a curve, the Narasimhan-Seshadri theorem
gives a cscK metric on the projectivisation of any polystable bundle (of
arbitrary rank) in any
K\"ahler class, and we are able to give an almost complete converse (there
is a small discrepancy for bundles which are strictly semistable but not
polystable until the
results of \cite{Do2, Do5,Mab} are improved to give K-polystability).

Other examples include unstable blow ups in Section \ref{sec:blowups}
and unstable rational manifolds in Section
\ref{sec:rational}.  In particular we give examples of
K\"ahler classes on surfaces with trivial automorphism group which do
not admit cscK metrics (\ref{rmk:folklore}).

One might hope that in the continuity method to find a cscK metric,
the multiplier ideal sheaf along which the $C^0$-estimates required
for closedness fail \cite{Na} defines a subscheme which slope
destabilises the variety. In particular if one could show this for canonically
polarised manifolds then Theorem \ref{thm:Kstableimpliesslopestable} combined
with Nadel's results would solve the K\"ahler-Einstein problem for Fano manifolds.

\subsection*{Notation and Terminology}\ \vskip 5pt

In this paper $(X,L)$ will be a smooth complex manifold of dimension
$n$ with a polarisation $L$ (\emph{i.e.}\ an ample line bundle on
$X$).  Furthermore $Z$ will denote an arbitrary subscheme of $X$
defined by an ideal sheaf $\I_Z$.  When $Z$ is smooth
$\nu_Z=(\I_Z/\I_Z^2)^*$ will denote its normal bundle.

The blow up along $Z$ is denoted by $\pi\colon\widehat X\to X$, with
exceptional divisor $E$.  Note that $\pi_* \O(-jE)=\I_{\!Z}^j$ for
$j\gg 0$.  For convenience we often suppress pullback maps and use the
same letter to denote a divisor and the associated line bundle.  For
example on $\widehat{X}$, we denote $(\pi^*L\otimes\O(-E))^{\otimes
  k}$ by $L^{\otimes k}(-kE)$.  The intersection product of divisors
$D_1,\ldots,D_n$ on $X$ is denoted by $\int _X c_1(D_1).\ldots
c_1(D_n)$, and this is abbreviated to $D_1.D_2\ldots D_n$ in sections
\ref{divisors} and \ref{sec:nonpositive}.


A $\mathbb Q$\,-divisor is a formal sum of divisors with rational
coefficients; some multiple is therefore a divisor with an associated
line bundle.  A $\mathbb Q$\,-divisor is said to be ample if it can be
written as a formal sum of ample divisors with positive rational
coefficients.  We recall that a nef line bundle (or divisor) is one
whose intersection with every curve in $X$ is nonnegative, and this
extends to $\mathbb Q$\,-line bundles.  By the Kleiman criterion
\cite{Kl} these divisors are precisely those in the closure of
the ample cone.  In notation like $H^0(L^{\otimes k})$ we always
tacitly restrict to those $k$ for which $L^{\otimes k}$ is an honest
line bundle.

Any finite-dimensional vector space $V$ with a $\C^\times$-action
splits into one-dimensional weight spaces $V=\bigoplus_iV_i$, where
$t\in\C^\times$ acts on $V_i$ by $t^{w_i}$. The integers $w_i$ are the
\emph{weights} of the action, and $w(V)=\sum_iw_i$ is the \emph{total
  weight} of the action; \emph{i.e.}\ the weight of the induced action
on the top exterior power $\Lambda^{\max}V$.
\subsection*{Acknowledgments}\ \vskip 5pt

We would most like to thank Simon Donaldson, who suggested looking at
stability of projective bundles $\PP(E)$ via the degeneration
mentioned in (\ref{ext1}). As we discuss there, the deformation to the
normal cone collapses to this degeneration in one special case, and it
was this case that motivated our study of the normal cone; without his
invaluable suggestion we would never have looked at stability of
varieties.  We also wish to thank S\'ebastien Boucksom, David
Calderbank, Joel Fine, Jun Li, Sean Paul, Gang Tian (both for useful
comments, and for inviting the first author to Princeton), Ben
Weinkove and Xiao-Wei Wang.

\section{Definition of K-stability}

Tian \cite{Ti2,Ti3} introduced a notion of K-stability using
differential geometry. Donaldson \cite{Do3} gave an algebro-geometric
definition that allows arbitrarily singular central fibre and which we
use here. The relation between the two is studied in \cite{PT}.

There is a strong formal link between K-stability and stability
notions in Geometric Invariant Theory; in particular the test
configurations defined below are what one gets by applying a one
parameter subgroup of projective linear transformations to the Kodaira
embedding of $(X,L^{\otimes k})$, and what we call the
Donaldson-Futaki invariant is the GIT weight of the induced action on
a certain line. We will not attempt to describe this further but
instead refer the interested reader to \cite{Do3, RT}.

\begin{defn} \label{def:testconfig}  Suppose that $(X,L)$ is a polarised variety with Hilbert polynomial $\P(k):=\chi(L^{\otimes k})$. A {\bf test configuration  with general fibre $(X,L)$} consists of
\begin{enumerate}
\item A flat projective family of $\mathbb Q$\,-polarised schemes $(\X,\L)\to
\mathbb C$.
\item An action of $\C^\times$ on $(\X,\L)$ covering the usual action
  of $\C^\times$ on $\C$,
  \end{enumerate}
  such that the fibre $(\X_t,\L|_{\X_t)})$ is
  isomorphic to $(X,L)$ for one, and so all, $t\in\C\backslash\{0\}$.
\end{defn}
The flatness condition is that the fibres $(X_t,\L_t)$ all
have the same Hilbert polynomial $\P(k)$ (\cite{Ha} Theorem III.9.9).
We call a test configuration a {\bf product configuration} if $\X\cong
X\times\C$, and a {\bf trivial configuration} if in addition
$\C^\times$ acts only on the second factor.  Since $0\in\mathbb C$ is
fixed, we get an induced action on the central fibre $(\X_0,\L|_{\X_0})$ and
hence on $H^0(\X_0,\L^{\otimes k}|_{\X_0})$ for all $k$.

\begin{defn}\label{def:futaki}
  Suppose $(\X,\L)$ is a test configuration with general fibre
  $(X,L)$.  Let $w(k)$ be the weight of the induced action on $ H^0(\X_0,\L|_{\X_0}^{\otimes k})$,
  which by the
  equivariant Riemann-Roch formula is a polynomial of degree
  $n+1$ for $k\gg 0$,  so there is an expansion
  $$
  \frac{w(k)}{k\P(k)} = f_0 + f_1k^{-1} + O(k^{-2}).$$
  We define the {\bf Donaldson-Futaki invariant} of a test configuration
  to be $F_1=-f_1$
  (so this has the opposite sign to the definition in \cite{Do3}).

  Writing $\P(k)=a_0k^n + a_1k^{n-1} + O(k^{n-2})$ and $w(k) =
  b_0k^{n+1} +b_1k^n + O(k^{n-1})$, then
\begin{equation}
 F_1 = \frac{b_0a_1-b_1a_0}{a_0^2}\,.\label{futaki_aibi}
\end{equation}
\end{defn}

\begin{defn} \label{defkstab}\
\begin{itemize}
\item We say that $(X,L)$ is {\bf algebraically K-stable} (resp.\ {\bf
    algebr\-aically K-semistable}) if for all non-trivial test
  configurations with general fibre $(X,L)$ we have $F_1>0$ (resp.\
  $F_1\ge 0$).
\item We say that $(X,L)$ is {\bf algebraically K-polystable} if it is K-semistable,
  and any test configuration with general fibre $(X,L)$ and $F_1=0$ is
  a product configuration. That is, the only instability arises from
  $\C^\times$-actions on $(X,L)$.
\end{itemize}
\end{defn}  

\begin{rmks}\label{rmk:kstable}\
  \begin{itemize}
  \item The property of being K-(semi/poly)stable is preserved under
    replacing $L$ by $L^{\otimes r}$, so makes sense when $L$ or $\L$ is an
    ample $\mathbb Q$\,-line bundle.  The definition of a test
    configuration given here differs from that in Definition 3.6 of
    \cite{RT}, but is the same after twisting $\L$ by some power.
  \item When the central fibre $(\X_0,\L|_{\X_0})$ is smooth, $F_1$ is, up
  to a constant, the usual Calabi-Futaki invariant, with respect to the class
  $c_1(L)$, of the vector field induced by the $S^1$-action \cite{Do3}.
  \item The Donaldson-Futaki invariant can be interpreted in terms of the
    Mumford weight function in Geometric Invariant Theory \cite{Do3}
    (see also Theorem 3.9 of \cite{RT}).
\end{itemize}
\end{rmks}


It is possible to strengthen the definition of K-stability.  We define
an {\bf analytic test configuration} with general fibre $(X,L)$
exactly the same way as we defined a test configuration, but allow
$\L$ to be an ample $\mathbb R$-divisor.  (By an ample $\mathbb
R$-divisor we mean a formal sum $\L=\sum_{i=1}^m \alpha_i D_i$ with
each $D_i$ an ample divisor and $\alpha_i$ a positive real; a $\mathbb
C^\times$-action on $\L$ is a choice of $\mathbb C^\times$-action on
each line bundle $\O(D_i)$ ).

For any test configuration the Donaldson-Futaki invariant can be
calculated using the equivariant Riemann-Roch theorem in terms of the equivariant
first Chern class of $\L$ with its $\C^\times$-action.
The resulting expression makes sense even if $\L$ is an ample $\mathbb
R$-divisor, and we take this to be the definition of $F_1$ in this case.


\begin{defn}\label{defan}\ \\
  We say that $(X,L)$ is {\bf  analytically K-stable} (resp.\ {\bf  analytically K-semi\-stable})
if for all non-trivial analytic test
  configurations $(\X,\L)$ with general fibre $(X,L)$ we have $F_1>
  0$ (resp.\ $F_1\ge 0$).   It is {\bf analytically K-polystable}
  if it is analytically K-semistable and any analytic test
  configuration with $F_1=0$ is a product configuration.
\end{defn}

\begin{rmk}
  As analytic K-semistability is equivalent to algebraic
  K-semistability, we will drop the qualifier when dealing with
  K-semi\-stability.
\end{rmk}

\subsection{Relationship to constant scalar curvature K\"ahler metrics}\label{sec:KstabandcscKmetrics}\ 
The precise conjecture relating K-stability to the
existence of cscK metrics is the following \cite{Y2, Ti3, Do3}.

\begin{conj}[Yau-Tian-Donaldson]\label{conj:cscK}
  Let $(X,L)$ be a polarised manifold.  Then there exists a constant
  scalar curvature K\"ahler metric in the class of $c_1(L)$ if and only
  if $(X,L)$ is K-polystable.
\end{conj}

One direction of this conjecture, that existence of a cscK metric
implies stability, has almost been proved: in \cite{Do5} it is shown
that a cscK metric implies K-semistability. Before that paper one had
a slightly weaker result by using balanced metrics: if $\aut(X)=0$
then the existence of a cscK metric implies that the Kodaira embedding
of $(X,L^{\otimes r})$ can be ``balanced" for $r\gg 0$ \cite{Do2}.
This implies it is asymptotically Chow stable \cite{Zh, P, Wa}, which
in turn implies that $(X,L)$ is  K-semistable (see for example Theorem
3.9 of \cite{RT}). Mabuchi \cite{Mab} extended this proof to manifolds
with non-discrete automorphism group satisfying a certain stability condition.

Another path to stability is through the K-energy, also called the
Mabuchi functional.  The existence of a cscK metric implies the
K-energy map is bounded from below \cite{Do4,CT} (resp. proper in the
K\"ahler-Einstein Fano case when $\aut(X)=0$ \cite{Ti3});
in turn this implies K-semistability \cite{PT} (resp.
K-stability).

We remark that a recent example in \cite{ACGT} suggests that algebraic
K-stability may not be enough to guarantee the existence of a
cscK metric, and that the stronger analytic definition of K-stability
may be required. (The authors wish to thank V.  Apostolov and
  D.\ Calderbank for discussions on this point).  Moreover it may be
that we have to allow non-projective central fibres (see Section
\ref{extensionkahler}).

It is expected that the deep results mentioned above proving stability
are not optimal, and that K-polystability can be proved.  However the
fact that a cscK metric implies K-semistability is enough to give a
new obstruction in terms of the subschemes of $X$ which we now
describe.

\section{Definition of slope stability}\label{sec:definitionofslope}
Fix a polarised manifold $(X,L)$ and write the Hilbert polynomial as
$$\P(k)=\chi(L^{\otimes k}) = a_0 k^n + a_1 k^{n-1} + O(k^{n-2}).$$

\begin{defn}
  The {\bf slope} of $(X,L)$ is $$\mu(X) = \mu(X,L) =
  \frac{a_1}{a_0}\,.$$
  By the Riemann-Roch theorem,
  $$a_0=\frac{1}{n!}\int_X c_1(L)^n,\quad\text{and}\quad a_1 =-\frac{1}{2(n-1)!}\int_X c_1(K_X).c_1(L)^{n-1},$$
so
  $$\mu(X) = -\frac{n\int_X c_1(K_X).c_1(L)^{n-1}}{2\int_X c_1(L)^n}\,.$$
\end{defn}

For a subscheme $Z$ of $X$ let $\widehat{X}$ be the blow up of $X$ along $Z$, with exceptional divisor $E$.

\begin{defn} \label{def:seshadri} The {\bf Seshadri constant} of $Z$ is
\begin{eqnarray*}
  \epsilon(Z)&=&\epsilon(Z,X,L) \\
&=& \sup\,\{ c: L^{\otimes k}\otimes\I_Z^{ck}\ \text{is globally generated for }k\gg0\}\\
&=&\sup\,\{ c : L(-cE)\,\,\text{ is ample on }\widehat X \}\\
&=&\max\,\{ c : L(-cE)\,\,\text{ is nef on }\widehat X \}\,.
\end{eqnarray*}

We say the global sections of $L\otimes \I_Z$ {\bf saturate} $\I_Z$ if they
generate the line bundle $L(-E)$ on $\widehat X$. This is weaker than
(\emph{i.e.}\ is implied by) $L\otimes \I_Z$ being globally generated (see \cite{RT} section 2).
\end{defn}

For fixed $x\in\mathbb Q$ define $a_i(x)$ by
\begin{equation} \label{aidef}
\chi(L^{\otimes k}(-xkE)) = a_0(x) k^n + a_1(x)k^{n-1} + O(k^{n-2}) \,\,\, k\gg
0,xk\in\mathbb N.
\end{equation}
Since $\chi(L^{\otimes k}(-rE))$ is a polynomial in $k$ and $r$ of total degree at
most $n$, $a_i(x)$ is a polynomial in $x$ of degree at most $n-i$,  so
can be extended to all of $\mathbb R$.  We have
\begin{equation}\label{definitionofa0x}
 a_0(x) = \frac{1}{n!}\int_{\widehat X} c_1(L(-xE))^n,
\end{equation}
and, when $Z$ is a codimension $p$ submanifold, by the Riemann-Roch formula
on
$\widehat{X}$,
\begin{equation}
a_1(x) = -\frac{1}{2(n-1)!}\int_{\widehat X} c_1(K_{\widehat X}).c_1(L(-xE))^{n-1},\label{definitionofa1x}
\end{equation}
where $K_{\widehat X}=K_X((p-1)E)$ is the canonical divisor of
$\widehat{X}$.

The $a_i(x)$ can also be defined in terms of the ideal sheaf of $Z$.  Fix
    $j_0$ such that $\pi_*(-jE)=\I_Z^j$ for all $j\ge j_0$ (when $Z$
    is smooth we can take $j_0=0$).  Then for $xk\in \mathbb N$,
    $x<\epsilon(Z)$ and $k\gg0$ (in particular $kx\ge j_0$), 
 \begin{equation} \label{alternativeHS}
 h^0(L^{\otimes k}\otimes\I_Z^{xk}) =\chi(L^{\otimes k}\otimes\I_Z^{xk})
    =a_0(x) k^n + a_1(x) k^{n-1} + O(k^{n-2}).
 \end{equation}
    Thus $a_0(0)=a_0$.  When $X$ and $Z$ are smooth,
    taking $j_0=0$ shows that we also have $a_1(0)=a_1$.  More
    generally this holds when $X$ is normal (\cite{RT} Remarks 4.21).

\begin{defn}
  The {\bf slope of $Z$ with respect to $c$} is 
  \begin{equation*}
  \mu_c(\I_Z) =\mu_c(\I_Z,L)= \frac{\int_0^c \left(a_1(x) + \frac{a_0'(x)}{2}
    \right)dx}{\int_0^c a_0(x) dx}\,.
\end{equation*}
\end{defn}

\begin{defn} \label{defslopestab}\ 
\begin{itemize}
\item We say that $(X,L)$ is {\bf slope semistable with respect to
    $Z$} if $\mu_c(\I_Z)\le\mu(X)$  for all   $c\in(0,\epsilon(Z)]$.
  
\item We say $(X,L)$ is {\bf slope stable with respect to $Z$} if
  $\mu_c(\I_Z,L)< \mu(X)$ for every $c\in(0,\epsilon(Z))$, \emph{and}
  for $c=\epsilon(Z)$ if $\epsilon(Z)$ is rational and the global
  sections $L^{\otimes k}\otimes\I_{\!Z}^{\epsilon(Z)k}$ saturate
  $\I_Z^{\epsilon(Z)k}$ for $k\gg0$.
  
\item We say $(X,L)$ is {\bf slope polystable with respect to $Z$} if
  it is slope semistable, and if $(Z,c)$ is any pair such that
  $\mu_c(\I_Z)=\mu(X)$, then $c=\epsilon(Z)\in\mathbb Q$ and, on the
  deformation to the normal cone (Section \ref{sec:normalcone}) of
  $Z$, $\L_c=L(-cP)$ is pulled back from a product test configuration
  $(X\times\C,L)$.
\item Finally $(X,L)$ is said to be {\bf slope (semi/poly)stable} if it is
  so with respect to all subschemes $Z$.
\end{itemize}
\end{defn}

\begin{rmk}
  The definition above of slope semistability agrees with that in
  \cite{RT}.  However the definitions given here of slope
  (poly)\-stability are slightly stronger as we require the relevant
  condition to hold even for irrational $c$.  Thus what we have
  defined as slope (poly)stability might more properly be referred to
  as \emph{analytic} slope (poly)stability, and is the notion relevant to
  the analytic K-stability of Definition \ref{defan}.
\end{rmk}

An example of a slope polystable variety is provided by $\PP^n$
(whose Fubini-Study metric is cscK). When $c=\epsilon(p),\ \mu_c(\I_p)=\mu(\PP^n)$
and the deformation to the normal cone (Section \ref{sec:normalcone}) of
a point $p\in\PP^n$ collapses to $\PP^n\times\C$, with a non-trivial
$\C^\times$-action with Donaldson-Futaki
invariant 0. Generalisations of this example are provided by the projective
bundles of (\ref{ext1}).

\begin{rmks}\label{rmk:onslopes} \
  \begin{itemize}
  \item We say that $Z$ {\bf destabilises} (resp.\ {\bf strictly
      destabilises}) if $(X,L)$ is not slope stable (resp.\ slope
    semistable) with respect to $Z$.
  \item Slope (semi/poly)stability is preserved under twisting $L$,
    since $\epsilon(Z,L^{\otimes r}) = r\epsilon(Z,L)$,
    $\mu(X,L)=r\mu(X,L^{\otimes r})$, and
    $\mu_{c}(\I_Z,L)=r\mu_{rc}(\I_Z,L^{\otimes r})$.
  \item If $0< x<\epsilon(Z)$ then from the fact that $L(-xE)$ is ample,
    \begin{eqnarray} \label{negative} \qquad\qquad
      a_0(x) &=& \frac{1}{n!}\int_{\widehat X}c_1(L(-xE))^n>0, \nonumber
      \\
      a_0'(x) &=& - \frac{1}{n!}\int_{\widehat X}c_1(L(-xE))^{n-1}.E<0.
    \end{eqnarray}
    In particular, for $0<c\le
    \epsilon(Z)$, $\int_0^c a_0(x)dx>0$ so $\mu_c(\I_Z)$ is finite.
  \item $\lim_{c\to0}\mu_c(\I_Z)=\frac{a_1(0)+a_0'(0)/2}{a_0}<\frac
    {a_1(0)}{a_0}$ by (\ref{negative}). For $X$ normal this is
    $\frac{a_1}{a_0}=\mu(X)$ (by Remarks 4.21 of \cite{RT}), so $(Z,c)$
    does not destabilise for small $c>0$. Therefore, on defining
    $\mu(\I_Z):=\max\_{0\le x\le c} \mu_c(\I_Z)$, slope semistability
    is equivalent to $\mu(\I_Z)\le\mu(X)$.  This is how it was
    presented in the Abstract.
    \end{itemize}
\end{rmks}
\begin{rmks}\
  \begin{itemize}
  \item In the slope inequality we may assume without loss of
    generality that $Z$ is not a thickening of any other subscheme.
    For if $Z=mZ'$, $m\ge 1$, then
    $\epsilon(Z) = \frac{1}{m}\epsilon(Z')$ and, as $a_0'(x)<0$
    (\ref{negative}),
\begin{eqnarray*} 
\mu_{c/m}(\I_{\!Z'}^m) &=& \mu_c(\I_{Z'}) + (m-1)\frac{\int_0^c a_0'(x)dx}
{2\int_0^c a_0(x)dx} < \mu_c(\I_{Z'}).
\end{eqnarray*}
\item If $Z$ strictly destabilises then so does one connected
  component of $Z$, and smooth points do not destabilise a smooth $X$
  (\cite{RT} Theorem 4.29).  Thus, for the strict inequality, we may
  assume without loss of generality that $Z$ is connected.
\end{itemize}
\end{rmks}

\begin{defn}\label{def:quotientslope}
  Let $\tilde{a}_i(x)$ be defined by
$$\chi(L^{\otimes k}\otimes\O_{xkZ})=\chi(L^{\otimes k}/(L^{\otimes k}\otimes\I^{xk}_{\!Z}))=\tilde{a}_0(x) k^n + \tilde{a}_1(x)k^{n-1} + O(k^{n-2}),$$ so
  $\tilde{a}_i(x) = a_i-a_i(x)$.  The {\bf quotient slope} of $Z$ with
  respect to $c$ is (in slightly misleading notation)
\begin{eqnarray} \label{qslope}
\mu_c(\O_Z)&=&\mu_c(\O_Z,L) \nonumber \\
&=&\frac{\int_0^c\left( \tilde{a}_1(x) + \frac{\tilde{a}'_0 (x)}{2}\right) dx}{\int_0^c \tilde{a}_0(x)dx}=\frac{\int_0^c\left(a_1(x) + \frac{a_0'(x)}{2} \right) dx - ca_1}{\int_0^c a_0(x)dx -ca_0}\,,
\end{eqnarray}
which is finite for $0<c\le \epsilon(Z)$. Notice that
$$
\mu_c(\I_Z)<\mu(X)\ \Longleftrightarrow\ \mu(X)<\mu_c(\O_Z)\ 
\Longleftrightarrow\ \mu_c(\I_Z)<\mu_c(\O_Z),
$$
due to the implications
$$
\frac{A}{B}< \frac{C}{D}\ \Longleftrightarrow\ 
\frac{C}{D}<\frac{C-A}{D-B}\ \Longleftrightarrow\ 
\frac{A}{B}<\frac{C-A}{D-B}
$$
for $0<B<D$, on setting $B=\int_0^c a_0(x) dx$ and $D=ca_0$ (so
$D-B=\int_0^c \tilde{a}_0(x)dx>0$).
\end{defn}

So slope stability can be phrased in terms of the quotient slope
$\mu_c(\O_Z)$. 

\begin{prop}\label{quotientslopeusingalpha}
  For fixed $x\in\mathbb Q_{\,>0}$, define $\alpha_i(x)$ by
\begin{equation}
  \chi(L^{\otimes k}\otimes \I^{xk}_Z /\I^{xk+1}_Z) =
\alpha_1(x)k^{n-1} + \alpha_2(x) k^{n-2} + O(k^{n-3})\label{alphai}
\end{equation}
for $k\gg 0, xk\in\mathbb N$.
(So if $Z$ is smooth with normal bundle $\nu_Z$ then
\begin{equation*}
\chi(L^{\otimes k}|_Z\otimes S^{xk}\nu_Z^*) =
\alpha_1(x)k^{n-1} + \alpha_2(x) k^{n-2} + O(k^{n-3}), \, k\gg 0,\ 
xk\in\mathbb N,
\end{equation*}
where $S^r(\cdot)$ denotes the $r$-th symmetric product.)

 Then
\begin{equation*}
 \mu_c(\O_Z) = \frac{\int_0^c (c-x)\alpha_2(x)dx + \frac c2\alpha_1(0)}
 {\int_0^c (c-x)\alpha_1(x) dx}\,.
\end{equation*}
\end{prop}
\begin{proof}
Fix $x>0$ and let $\bar{x}=x+1/k$.   Clearly
 \begin{eqnarray*}
\chi(L^{\otimes k} \otimes \I_Z^{xk} /\I_Z^{xk+1}) &=& \chi(L^{\otimes k} \otimes \I_Z^{xk}) - \chi(L^{\otimes k} \otimes \I_Z^{xk+1}). 
\end{eqnarray*}
By (\ref{alternativeHS}) and the Taylor expansion of $a_i(x)$ this
equals, for $k\gg 0$,
\begin{align*}
[a_0(x) &- a_0(\bar{x})]k^n + [a_1(x) - a_1(\bar{x})]k^{n-1}  +\cdots\\
&= -a_0'(x)k^{n-1} - \frac{a_0''(x)}{2}k^{n-2} -a_1'(x) k^{n-2} + O(k^{n-3}).
 \end{align*}
 (Note that this holds when $n=1$ for then $a_0''(x)=a_1'(x)=0$).  Hence
 \begin{equation}
\alpha_1(x) = -a_0'(x) \quad \text{and}\quad \alpha_2(x) =-a_1'(x) -
 \frac{a_0''(x)}{2}\,.\label{eq:alphasandas}
\end{equation}
  Thus the denominator of the quotient slope of
 $Z$ is $$\int_0^c \tilde{a}_0(x)dx = \int_0^c \int_0^x \alpha_1(y) dy
 dx = \int_0^c (c-x)\alpha_1(x) dx.$$  The calculation for the
 numerator is similar.
\end{proof}

\section{Slope stability as a necessary condition for K-stability}\label{sec:kstabilityimpliesslopestability}

\subsection{Deformation to the normal cone}
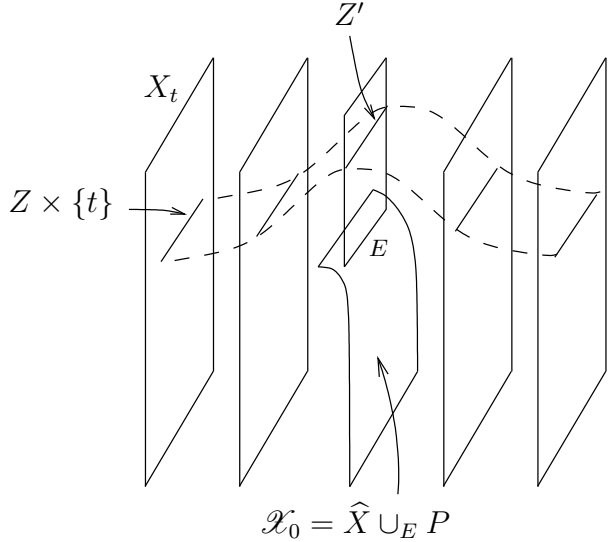
\begin{figure}[h]
  \center{\input{cone.pstex_t}}
  \caption{The deformation to the normal cone of $Z$. \label{fig:normalcone}}
\end{figure}

\ \vskip 5pt \label{sec:normalcone}
Fix a subscheme $Z$ of $(X,L)$.  Let $\X$ be the deformation to the
normal cone of $Z$, so $\X$ is the blow up of $X\times\mathbb C$ along
$Z\times\{0\}$, and denote the exceptional divisor by $P$.  The
central fibre $\X_0$ is isomorphic to the blow up $\widehat X$ glued to
$P$ along $E$ (see Figure \ref{fig:normalcone}).  When $Z$ is a
submanifold, $E=\PP(\nu)$ and $P$ is isomorphic to the projective
completion of the normal bundle of $Z$, \emph{i.e.}\ 
$P=\PP(\nu_Z\oplus\underline{\mathbb C})$, with a copy $Z':=
\PP(\underline{\mathbb C})$ of $Z$ as its zero section.

Consider the product action on $(X\times \mathbb C,L)$ (where as usual we
suppress the pullback map on $L$), which acts trivially
on $(X,L)$ but scales $\mathbb C$ with weight 1.  This fixes $Z\times
\{0\}$ and so induces an action on $\X$ and on $P$.  The induced
action on the central fibre $\X_0=\widehat X \cup_{E} P$ is trivial on
$\widehat X$, and $\lambda\in\mathbb C^\times$ acts on
$P=\PP(\nu_Z\oplus\underline{\mathbb C})$ as $\text{diag}(1,\lambda)$.

We define a $\mathbb Q$\,-line bundle on $\X$ by $\L_c=L(-cP)$ for $c\in\mathbb
Q$.

\begin{lem}\label{thm:degentonormalconeisatestconfig}
  For rational $c\in(0,\epsilon(Z))$, the line bundle $\L_c$ is ample.
\end{lem}
\begin{proof}
  Let $p\colon\X\to X$ be the composition of the projections, and let
  $c=r/q$. Choose $q$ and $r$ sufficiently large so that $L^{\otimes
    q}$ and $L^{\otimes q}(-rE)$ are globally generated.  Then away from
  $Z'=\PP(\underline{\mathbb C})\subset P$, the line bundle
  $\L_c^{\otimes q}=L^{\otimes q}(-rP)$ is generated by $p^*
  H^0(L^{\otimes q}\otimes\I_Z^r)$, while on $Z'$, it is generated by
  $t^r p^*H^0(L^{\otimes q})$.  That is, $\L_c^{\otimes q}$ is globally
  generated for all $0<c<\epsilon(Z)$ and so nef for
  $c\in[0,\epsilon(Z)]$.  Since $\L_c$ is ample for $c$ sufficiently
  small, the fact that the ample cone is convex and the interior of
  the nef cone \cite{Kl} implies that $\L_c$ is ample for rational
  $0<c<\epsilon(Z)$.
\end{proof}

\subsection{Slope stability as an obstruction to K-stability}\label{sec:slopestabilityasanobstruction}\
\nopagebreak \vskip 5pt 

Slope stability with respect to $Z$ is precisely K-stability
restricted to test configurations arising from the degeneration to the
normal cone of $Z$.  In (\cite{RT} Theorem 4.18) it is shown that
K-semistability implies slope semistability.  Moreover, using the
algebraic definitions of slope (polystability) in \cite{RT} it is also
shown that K-(poly)stability implies slope (poly)stability.  Here we
give a proof of the part of this result which is sufficient for our
examples and applications to cscK metrics.

\begin{thm}\label{thm:Kstableimpliesslopestable}
  Suppose $(X,L)$ is K-semistable. Then it is slope semistable with
  respect to any smooth subscheme $Z$.
\end{thm}
\begin{proof}
  We need to show that $\mu_c(\O_Z)\ge \mu(X)$ for all
  $0<c\le\epsilon(Z)$.  By continuity of $\mu$ with respect to $c$ it
  is sufficient to consider rational $c<\epsilon(Z)$.  So $\L_c$ is
  ample (\ref{thm:degentonormalconeisatestconfig}) and hence so is
  $\L_c|_{\X_0}$.
  
  By the definition of the blow up in $\I_{Z\times\{0\}}\subset
  \O_{X\times\C}$ (\emph{i.e.}\ $\I_Z+(t)\subset\mathbb C[t]\otimes
  \mathcal O_X$, where
  $t$ is the coordinate on $\mathbb C$), for $k\gg 0$ and $ck\in\mathbb N$,
  \begin{eqnarray} \label{fns}
  H^0(\X,\L_c^{\otimes k}) &=& H^0(\X,(L(-cE))^{\otimes k}) \nonumber \\
&=&H^0(X\times\C,L^{\otimes k}\otimes
  \I_{Z\times\{0\}}^{ck})   \nonumber \\
   &=&H^0(X\times\C,L^{\otimes k}\otimes(\I_Z+(t))^{ck})\\
&=&\bigoplus_{i=1}^{ck}
  t^{ck-i}H^0(X,L^{\otimes k}\otimes\I_Z^i)\ \oplus\,t^{ck} \mathbb C[t] H^0(L^{\otimes k}). \nonumber
  \end{eqnarray}
Similarly, for $k$ sufficiently large and $j\ge1$,
$$0=H^j(\X,\L_c^{\otimes k})=\bigoplus_{i=1}^{ck}t^{ck-i}H^j(X,L^{\otimes k}\otimes\I_Z^i) \oplus\,t^{ck}
\mathbb C[t] H^j(L^{\otimes k}),$$ so that
\begin{equation}
H^j(L^{\otimes k}\otimes\I_Z^i)=0 \quad \text{for}\,\, j\ge1,\ k\gg 0,\ ck\in\mathbb N,\ i=0,\ldots,ck. \label{vanishing}
\end{equation}
Thus, for instance, $H^0(L^{\otimes k}\otimes\I_Z^i)\big/
  H^0(L^{\otimes k}\otimes\I_Z^{i+1})=H^0\big(L^{\otimes k}\otimes\big(\I_Z^i/\I_Z^{i+1}\big)\big)
  =H^0(L^{\otimes k}|_Z\otimes S^i\nu_Z^*)$. So from (\ref{fns}) we get the splitting,
for $k\gg0$,
\begin{multline}
H^0(\L_c^{\otimes k}|_{\X_0}) = H^0(\X,\L_c^{\otimes k})\big/tH^0(\X,\L_c^{\otimes k})= \\
H^0(X,L^{\otimes k}\otimes\I_Z^{ck}) \oplus \bigoplus_{i=0}^{ck-1}
t^{ck-i} H^0(L^{\otimes k}|_Z\otimes S^{i}\nu_Z^*), \label{eq:spaceofsections}
\end{multline}
of the functions on $\X_0$ into those on $\widehat X$ and the polynomials
on the $\nu_Z$-fibres of $P$.
In particular $h^0(\L_c^{\otimes k}|_{\X_0})$ equals
\begin{align*}
  h^0(L^{\otimes k}\otimes\I_Z^{ck}) &+
\sum_{i=0}^{ck-1}\left(h^0(L^{\otimes k}\otimes\I_Z^{i})
  -h^0(L^{\otimes k}\otimes\I_Z^{i+1})\right)\\
&= h^0(L^{\otimes k})=\P(k).
\end{align*}
This proves flatness, so $(\X,\L_c)$ is a test configuration
with general fibre $(X,L)$.

Now $\mathbb C^\times$ acts
trivially on $(X,L)$ and so also on $\nu^*_Z$ and $L|_Z$, but with weight
$-1$ on $t$, so (\ref{eq:spaceofsections}) is also the weight space decomposition
of $H^0(\L_c^{\otimes k}|_{\X_0})$ into the pieces of weight $-(ck-i)$. Thus the total
weight of the action on $H^0(\L_c^{\otimes k}|_{\X_0})$ is
\begin{eqnarray*}
w(k) &=& - \sum_{i=0}^{ck-1} (ck-i) h^0(L^{\otimes k}|_Z\otimes S^i \nu_Z^*)\\
&=& -\sum_{i=0}^{ck-1} (ck-i) \chi(L^{\otimes k}|_Z\otimes S^i \nu_Z^*)\\
&=& -\sum_{i=0}^{ck-1} (ck-i)\left(\alpha_1(i/k) k^{n-1} + \alpha_2(i/k)
k^{n-2} + O(k^{n-3})\right),
\end{eqnarray*}
using the fact that $H^j(L^{\otimes k}\otimes S^i \nu_Z^*)=0$ for $j>0,\ k\gg 0,\ ck\in\mathbb
N,\ i=0,\ldots,ck-1$ (by (\ref{vanishing}) and $S^i \nu_Z^*=\I_Z^i/\I_Z^{i+1}$).
Here the $\alpha_i$ are as in (\ref{alphai}) and (\ref{eq:alphasandas}).
The $k^{n+1}$ and
$k^n$ terms of $w(k)$ can be calculated using the trapezium rule (Lemma
\ref{trapezium}), giving $w(k) = b_0k^{n+1} + b_1 k^n + O(k^{n-1})$,
where
\begin{eqnarray}
  \label{eq:bi}
  b_0 &=& -\int_0^c (c-x) \alpha_1(x) dx=\int_0^c\!a_0(x)dx-ca_0, \nonumber \\
b_1 &=& -\int_0^c \left((c-x) \alpha_2(x) + \frac{\alpha_1(0)}{2}\right)dx\\
&=&\int_0^c\!\left(a_1(x)+\frac{a_0'(x)}2\right)dx-ca_1,\nonumber
\end{eqnarray}
where each line follows from integration by parts and (\ref{eq:alphasandas}).

As $(X,L)$ is assumed to be K-semistable, the Donaldson-Futaki
invariant $F_1$ of the test configuration $(\X,\L_c)$ is nonnegative so
$$
0\le F_1=\frac{1}{a_0^2}(b_0a_1-b_1a_0) =
\frac{-b_0}{a_0}\left(\frac{b_1}{b_0}-\frac{a_1}{a_0}\right)=
\frac{-b_0}{a_0}\left(\mu_c(\O_Z)-\mu(X)\right),$$ where the last
equality uses (\ref{def:quotientslope}) and (\ref{eq:bi}).  Using
(\ref{eq:alphasandas}, \ref{rmk:onslopes}), $\alpha_1(x) = -a_0'(x)$ is
positive for $0<x<\epsilon(Z)$.  By equation
(\ref{eq:bi}) this shows that $b_0<0$, and hence $\mu_c(\O_Z)<\mu(X)$ as
required.
\end{proof}
\begin{lem}\label{trapezium}
  Let $f(x)$ be a polynomial.  Then
  $$\sum_{i=0}^{ck-1} (ck-i)f(i/k) = \int_0^c\left( k^2(c-x)f(x) +
    \frac{k}{2} f(0)\right) dx + O(k^0).$$
\end{lem}
\begin{proof}

  If $f(x)=\alpha$ is constant then both sides
  equal $\frac{\alpha}{2}ck(ck+1)$. So by linearity we may assume
  $f(x)=x^m,\ m\ge 1$.  Using $\sum_{i=0}^k i^m = \frac{1}{m+1}k^{m+1}
  + \frac{1}{2} k^{m} + O(k^{m-1})$ we get
  \begin{eqnarray*}
\sum_{i=0}^{ck-1} (ck-i)f(i/k) &=& k^{-m}\sum_{i=0}^{ck} (ck-i)i^m \\
&=& \int_0^c k^2(c-x)x^m dx + O(k^0),
\end{eqnarray*}
as required.
\end{proof}

Although we will not use it, we indicate how this result extends to
K-stability.

\begin{thm}
  Suppose $(X,L)$ is analytically K-stable.  Then it is slope stable
  with respect to any smooth subscheme.
\end{thm}
\begin{proof}
  Suppose that $\mu_c(\O_Z) = \mu(X)$ for some $0<c<\epsilon(Z)$ (with
  $c$ possibly irrational).  By convexity of the ample cone and
  \eqref{thm:degentonormalconeisatestconfig} $\L_c$ is ample and thus
  the degeneration to the normal cone $(\X,\L_c)$ is an analytic test
  configuration.  For rational $d$ close to $c$ we have from the
  previous proof a test configuration $(\X,\L_d)$ with Futaki
  invariant
$$F_1(d) = \frac{-b_0}{a_0}\left(\mu_{d}(\O_Z)-\mu(X)\right).$$  Thus the Futaki invariant of $(\X,\L_c)$ is
$$F_1 = \lim_{d\to c} F_1(d) =\lim_{d\to c}
\frac{-b_0}{a_0}\left(\mu_{d}(\O_Z)-\mu(X)\right)=0$$
since
$\lim_{d\to c} \mu_{d}(\O_Z)-\mu(X) = \mu_c(\O_Z)-\mu(X)=0$.  Thus
$(X,L)$ is not analytically K-stable.

Now if $c=\epsilon(X)$ is
rational and $\mu_c(\O_Z)=\mu(X)$ then it is shown in \cite{RT}
Theorem 4.18 that $(X,L)$ is not algebraically K-stable, and thus is
not analytically K-stable either.
\end{proof}

\subsection{Toric test configurations}\ \vskip 5pt

For toric varieties we can relate Donaldson's weight computation \cite{Do3}
to ours by an application of Fubini's theorem; \emph{i.e.}\ a change of order
of integration.
  Let $X_P=(X,L)$ be toric, defined by an integral polytope
  $P\subset\R^n$ such that $kP\cap\mathbb Z^n\cong H^0(X,L^{\otimes k})$.  Let
  $f\colon P\to\R$ be a strictly positive, rational, concave and piecewise
  linear function. Then the polytope
  $$
  Q=\{(p,t)\in P\times \mathbb R :  0\le t\le f(p)\},
  $$
  defines a toric variety with a $\Q$\,-polarisation $\L$, a $\C^\times$-action
  and an equivariant flat map to $\PP^1$. Removing the fibre over $\{\infty\}\in\PP^1$
  gives a test configuration $(\X,\L)$ with general
  fibre $(X,L)$ and $\C^\times$ acting on the section $(s,i)\in kQ\cap
  \mathbb Z^{n+1}\cong H^0(\X,\L^{\otimes k})$ with weight $-i$.
  
  Let $\#(kQ)$ denote the number of lattice points in $kQ$.  When $f$
  is \emph{integral}, Donaldson \cite{Do3} shows that the weight of this
  degeneration is
  $w_k=\#(kP)-\#(kQ)=b_0 k^{n+1} + b_1 k^n + O(k^{n-1})$ where
  \begin{equation}
  b_0 = -\int_P f d\mu=-\vol(Q), \qquad b_1=-\frac{1}{2}\int_{\partial
  P} f d\sigma. \label{donaldsonscalculationoftoricweight}
\end{equation}
Here $d\mu$ is the standard measure on $\mathbb R^n$ and
$d\sigma$ is defined by requiring that on any face of $P$ given by a
primitive integral conormal vector $h:\R^n\to \R$, we have $d\sigma\wedge
dh=\pm d\mu$.
 
Any toric subvariety of $X$ is defined by a face of $P$.  Such a face is
an intersection of codimension 1 faces.  Pick primitive integral
conormal vectors $\{f_i\}_{i=0}^m$ to the faces, with their signs chosen
so that $f_i\ge 0$ on $P$. Then the ideal of the subvariety is generated by the monomials
$$\big\{p\in P\cap\mathbb Z^n\colon f_i(p)\ge 1 \text{ for some } i\big\}
=\left\{p\in P\cap\mathbb Z^n\colon
\sum_{i=1}^m f_i(p)\ge 1\right\},$$
since $f_i(p)\ge 0$ for all $p\in P$. Therefore, more generally, the ideal
of any integrally closed toric \emph{subscheme} $Z$ (with multiplicities
$m_i$ in the direction of $f_i$) is generated by the monomials
\begin{equation} \label{toricideal}
\left\{p\in P\cap\mathbb Z^n\colon\sum_{i=1}^m\frac{f_i(p)}{m_i}\ge 1\right\}.
\end{equation}
(We have lost nothing by passing to the integral closure of $\I_Z$; this
corresponds to taking the
normalisation of the deformation to the normal cone of $Z$, and in testing
K- or slope stability one need only consider normal test configurations
(by Proposition 5.1 of \cite{RT}) since their Futaki invariants are smaller and so less stable.)

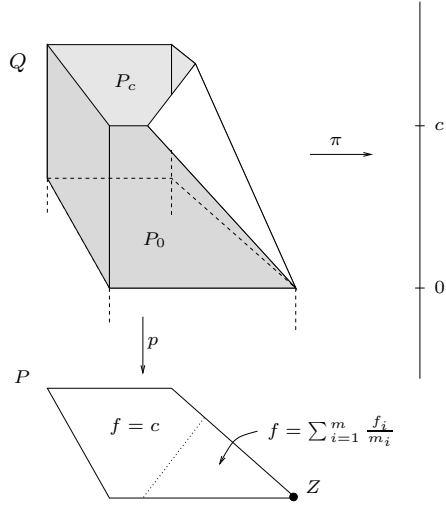
\begin{figure}[h]
  \center{\input{toriccone.pstex_t}}
  \caption{Toric representation $Q=\,$graph$(f)$ of the deformation to the
normal cone of $Z\subset X_P$. \label{toriccone}}
\end{figure}

The deformation to the normal cone $(\X,\L_c)$ of this $Z$ corresponds to
taking the positive, rational, concave, piecewise linear function $f=
\min\big(c,\sum_{i=1}^m\frac{f_i}{m_i}\big)$ in Donaldson's construction
(see Figure \ref{toriccone}, which sho\-uld of course be compared to Figure
\ref{fig:normalcone}).
(This $f$ is $\ge0$ but not everywhere $>0$, so to get the right geometry
we must add a positive constant to it. Since the resulting Donaldson-Futaki
invariant is independent
of the constant, we calculate without it.) Thus $f\colon P\to[0,c]$ and
\begin{equation}
f^{-1}[x,c]=P_x:=\left\{p\in P\colon\sum_{i=1}^m\frac{f_i(p)}{m_i}\ge
x\right\},
\end{equation}
which from (\ref{toricideal}) is seen to have integral points in $\frac1k\mathbb
Z^n$ which form a basis for $H^0(X,L^{\otimes k}\otimes\I_Z^{xk})$. So comparing coefficients
in $h^0(X,L^{\otimes k}\otimes\I_Z^{xk})=\#(kP_x)=\vol(P_x)k^n+\frac12
\vol(\partial P_x)k^{n-1}+O(k^{n-2})$ yields
\begin{equation} \label{toricformulae}
\vol(P_x)=a_0(x),\qquad \vol(\partial P_x)=2a_1(x).
\end{equation}

To relate Donaldson's weight formula (\ref{donaldsonscalculationoftoricweight})
to ours (\ref{eq:bi}) we change
the order of integration with respect to the two projections $Q\rt{p}P$ and
$Q\rt{\pi}\R$. That is, using (\ref{toricformulae}),
\begin{equation} \label{b0}
\int_Pf=\int_Pp_*1=\int_Q1=\int_0^c\pi_*1=
\int_0^c\vol(P_x)dx=\int_0^ca_0(x)dx.
\end{equation}

Similarly we can compute the volume of $\partial Q\backslash(P_0\cup P_c)$
(the ``sides" of $Q$) as $\int_0^c(\pi|\_{\partial Q})_*1=\int_0^c\vol(\partial
P_x)dx=2\int_0^ca_1(x)dx$ using (\ref{toricformulae}). Now $\vol(P_0)=a_0$
and $\vol(P_c)=a_0(c)$, so $\vol(\partial Q)=2\int_0^ca_1(x)dx+a_0+a_0(c)$.
But this can be computed differently as the volume of the shaded area in
Figure \ref{toriccone}, plus the volumes of the top and bottom.
If $f$ is integral (which is always assumed in \cite{Do3}) then the top and
bottom are (piecewise) integrally affine isomorphic
by the projection $p$; equivalently they have the same number of integral
points and so the same volume $a_0$. If $f$ is rational there are less
points on the top face, so the result has larger Futaki invariant, \emph{i.e.}\
it is more
stable. Alternatively we can multiply $f$ by an integer $N$ to make it integral;
this corresponds to taking the $N$th power of the $\C^\times$-action and
normalising the resulting test configuration.
Again, by Proposition 5.1 of \cite{RT}, this has Futaki invariant more
stable than ($N$ times) the old Futaki invariant. So either way we may as
well assume, like Donaldson, that $f$ is integral.

The area of the shaded region is computed
by $\int_{\partial P}f$, so we have found that
$$
2\int_0^ca_1(x)dx+a_0+a_0(c)=\int_{\partial P}f+2a_0,
$$
and so
\begin{equation} \label{b1}
\frac12\int_{\partial P}f=\int_0^c\!a_1(x)dx+\frac12(a_0(c)-a_0)=
\int_0^c\!\left(a_1(x)+\frac{a_0'(x)}2\right)dx.
\end{equation}
(\ref{b0}) and (\ref{b1}) differ from $-b_0$ and $-b_1$ in (\ref{eq:bi})
by
$ca_0$ and $ca_1$ respectively, which cancel in the Futaki invariant (\ref{futaki_aibi})
(or could be removed by adding $c$ to $f$). So we recover Donaldson's formulae
in this case.

\subsection{Extension to K\"ahler manifolds}\label{extensionkahler}\ \vskip 5pt

The definition of K-(poly/semi)stability given in
\eqref{defkstab} cannot be defined when the K\"ahler class is not
rational, but slope (poly/semi)\-stability can.  The same issue arises
for bundles; GIT cannot construct moduli for bundles over
non-projective manifolds, but the slope criterion for stability
generalises to all K\"ahler manifolds and Uhlenbeck-Yau proved that it
is equivalent to the existence of a HYM connection in this generality.

To define slope stability we must define the slope of an analytic
subspace $Z$ of a K\"ahler manifold $(X,\omega)$. We work on the blow
up $\pi\colon\widehat X\to X$ of $X$ in $Z$, with exceptional divisor
$E$. By the singular Hirzebruch-Riemann-Roch formula for analytic
spaces \cite{Ful} we can define a Todd \emph{homology} class of
$\widehat X$, and then define the polynomials $a_i(x)$ by the formula
$$
\int_{\mathrm{Td}\,(\widehat X)}\exp(k\pi^*\omega-xke)
=a_0(x)k^n+a_1(x)k^{n-1}+O(k^{n-2}).
$$
Here $e$ denotes any differential form Poincar\'e dual to $E$, and
when $X$ is projective with $\omega=c_1(L)$ this gives the same definition
as (\ref{aidef}). In particular $a_0(x)=\frac1{n!}\int_{\widehat X}(\pi^*\omega-xe)^n$,
while we can write $a_1(x)$ in terms of any resolution of singularities
$p\colon\,\overline{\!X}\to\widehat X$:
$$
a_1(x)=\frac1{2(n-1)!}\int_{\overline{\!X}}((p\comp\pi)^*\omega-xp^*e)^{n-1}
c_1(\overline{\!X}).
$$
Take any $c>0$ such that $\omega-ce$ has nonnegative volume
on any analytic subvariety of $\widehat X$ (if $\widehat X$ is smooth then
this is the condition that $\omega-ce$ be in the closure of the K\"ahler
cone of $\widehat X$ \cite{DP}). Then define the slope of $\I_Z$, with respect
to $\omega$ and $c$, as before:
$$
\mu_c(\I_Z):=\frac{\int_0^c\left(a_1(x)+\frac{a_0'(x)}2\right)dx}
{\int_0^ca_0(x)dx}\,.
$$
We say that $X$ is slope semistable if $\mu_c(\I_Z)\le \mu(X)$ for
all proper $Z\subset X$ and $c$ such that $\omega-ce$ has nonnegative
volume on any analytic subvariety of $\widehat X$.  For
slope stability we require that $\mu_c(\I_Z)< \mu(X)$ for all $c$ such
that $\omega-ce$ is the pullback of a K\"ahler form on a K\"ahler
variety.  We define $X$ to be slope polystable if it is slope
semistable and $\mu_c(\I_Z,\omega)=\mu(X)$ implies that $\omega-ce$ on
the deformation to the normal cone of $Z$ is pulled back from a map to the
product $X\times\C$.

Since the $\mathbb C^\times$-action on the degeneration to the normal
cone is trivial on the central fibre except on the component $P$, one
can use the localisation formula on $P$ to calculate the Calabi-Futaki
invariant in terms of the resulting vector field.  This gives an
alternate, but fundamentally equivalent, definition of slope for an
analytic $Z$ in a K\"ahler manifold. Then one would expect that the
usual argument that the derivative of the Mabuchi functional is the
Futaki invariant of the central fibre (defined in terms of the vector
field) should show that if $X$ is not slope polystable then the
Mabuchi functional is not proper, and so the class $[\omega]$ does not
admit a cscK metric.

\section{Examples}

\subsection{Slope of divisors and curves} \label{divisors}\ \vskip 5pt
\noindent
\begin{thm} \label{constants}
  Let $(X,L)$ be a polarised manifold of dimension $n\ge 2$ and suppose
  that $Z$ is a smooth curve in $X$ of genus $g$ with normal bundle
  $\nu_Z$. Then
  $$
  \mu_c(\O_Z)=\frac{n^2(n^2-1)(L.Z) - cn(n+1)[(n-2)c_1(\nu_Z) +
    2(g-1)]}{2nc[(n+1)(L.Z) - cc_1(\nu_Z)]}\,.
  $$
\end{thm}
 \begin{proof}
   The Riemann-Roch theorem for curves yields
    \begin{eqnarray*}
    \chi(L^{\otimes k}|_Z\otimes S^{xk} \nu_Z^*) &=& \rk S^{xk}
    \nu_Z \cdot   (kL.Z - \frac{xk c_1(\nu_Z)}{n-1}
 + 1-g) \\
&=& \alpha_1(x)k^{n-1} + \alpha_2(x)k^{n-2} + O(k^{n-3}),
  \end{eqnarray*}
  since $\frac{c_1(S^i\nu^*_Z)}{\rk S^i \nu_Z}= - i
  \frac{c_1(\nu_Z)}{n-1} $.  Now $\rk S^{xk}\nu_Z=\binom{xk+n-2}{n-2}$ equals
  $$\frac{1}{(n-2)!} \left(
    x^{n-2}k^{n-2} + \frac{(n-2)(n-1)}{2} x^{n-3}k^{n-3} +
    O(k^{n-4})\right).$$
  (This makes sense even if $n=2$ as in that case the $k^{n-3}$ term vanishes.) Thus
\begin{eqnarray*}
  \alpha_1(x) &=& \frac{x^{n-2}}{(n-2)!}\left(L.Z-\frac{xc_1(\nu_Z)}{n-1}
\right),\\
  \alpha_2(x) &=& \frac{x^{n-3}}{(n-2)!} \left(\frac{(n-2)(n-1)}{2}\left(L.Z-\frac{xc_1(\nu_Z)}{n-1}\right)+x(1-g)\right).
\end{eqnarray*}
Integration and rearranging (\ref{quotientslopeusingalpha}) gives the
formula for $\mu_c(\O_Z)$.
\end{proof}

\begin{thm} \label{divisor}
  Suppose that $Z$ is a divisor in $(X,L)$.  Then
  $$
  \mu_c(\O_Z)
  =\frac{n\left(L^{n-1}.Z-\sum_{j=1}^{n-1}\binom{n-1}{j}\frac{(-c)^j}{j+1}
      L^{n-1-j}.Z^j.(K_X(Z))\right)}{2\sum_{j=1}^n
    \binom{n}{j}\frac{(-c)^j}{j+1} L^{n-j}.Z^j}\,.
  $$
\end{thm}

\begin{proof}
  As $Z$ is a divisor, $\widehat X=X$ so
  (\ref{definitionofa0x}, \ref{definitionofa1x})
\begin{equation*}
  \tilde{a}_0(x)=a_0-a_0(x)=\frac{1}{n!}( L^n - (L-xZ)^n)= -\frac{1}{n!} \sum_{j=1}^n \binom{n}{j} L^{n-j}. (-xZ)^j,
\end{equation*}
and
\begin{eqnarray*}
  \tilde{a}_1(x) + \frac{\tilde{a}_0'(x)}{2} &=& \frac{1}{2(n-1)!} \left
  (-K_X. L^{n-1} + (K_X(Z)).(L-xZ)^{n-1} \right) \\
&=& \frac{1}{2(n-1)!} \sum_{j=1}^{n-1} \binom{n-1}{j}
  L^{n-1-j}. (-xZ)^j . (K_X(Z)) \\
&&+\, \frac{1}{2(n-1)!} L^{n-1}. Z.
\end{eqnarray*}
Integrating these expressions gives the required formula.
\end{proof}

The formulae (\ref{constants}) and (\ref{divisor}) agree for curves in
surfaces; the result simplifies to the following.

\begin{cor} \label{curvsurf}
  Let $Z$ be a smooth curve in a smooth polarised surface $(X,L)$.
  Then
\begin{eqnarray*}
\mu(X)&=&-\frac{K_X.L}{L^2}\,,\\
\mu_c(\O_Z)&=&\frac{3[2L.Z - c(K_X.Z + Z^2)]}{2c(3L.Z-c Z^2)}\,.
\end{eqnarray*}
If $Z$ is a smooth rational curve then
\begin{eqnarray*}
\mu_c(\O_Z)=\frac{3(L.Z +c)}{c(3L.Z-c Z^2)}\,.
\end{eqnarray*}
\end{cor}

We use these formulae in Section \ref{sec:rational} to give examples of
unstable rational surfaces.

\subsection{Manifolds with nonpositive first Chern class}\label{sec:nonpositive}\ \vskip 5pt
\noindent 
The existence of K\"ahler-Einstein metrics when $c_1(X)\le 0$ gives
K-semistability in these cases by the results of Donaldson \cite{Do5}.
We give a direct proof that smooth subschemes $Z\subset X$ do not
destabilise; a more general proof for arbitrary $Z$ (and extended to
varieties $X$ with canonical singularities) is given in Theorem 8.4 of
\cite{RT}.

\begin{thm} \label{CY} A polarised manifold $(X,L)$ is slope stable with respect to smooth subschemes if either
  \begin{enumerate}  
  \item $K_X$ is numerically trivial, or
  \item $K_X$ is ample and $L$ is a multiple of $K_X$.
\end{enumerate}
\end{thm}

\begin{proof}
In both cases, $K_X\sim\alpha
L$ is numerically equivalent to a
  nonnegative multiple $\alpha\ge0$ of the polarisation. So
  $\mu(X)=a_1/a_0=-nK_X.L^{n-1}/2L^n=-n\alpha/2$.
  If $Z\subset X$ is a codimension $p$ submanifold, 
  the canonical
  divisor of the blow up is $K_{\widehat X}=K_X((p-1)E)$.  Letting $L_x:=L(-xE)$,
\begin{eqnarray*}
-\mu(X)a_0(x)+a_1(x)&=&\frac{\alpha}{2(n-1)!}L_x^n-\frac1{2(n-1)!}
K_{\widehat X}.L_x^{n-1}\\
&=&-\frac1{2(n-1)!}\,(\alpha x+p-1)L_x^{n-1}.E\ \le\ 0,
\end{eqnarray*}
since $L_x$ is nef for $x\in(0,\epsilon(Z))$. As $a_0'(x)<0$ (\ref{negative}),
integration gives
$$
-\mu(X)\int_0^ca_0(x)dx+\int_0^ca_1(x)+\frac{a_0'(x)}2dx<0
\quad\text{for}\ c\in(0,\epsilon(Z)].
$$
Rearranging this gives slope stability, $\mu_c(\I_Z)<\mu(X)$.
\end{proof}

With more work these results can be extended to show that when $K_X$ is nef and $K_X^n>0$, then $X$ is slope stable for $L$ sufficiently close to $K_X$. More precisely, using additive notation ($aL+bK:=L^{\otimes a}\otimes K^{\otimes
b}$) for line bundles,

\begin{thm} \label{canonicalbundleimpliesstability}
  Fix a polarised manifold $(X,L)$ with $K_X$ nef and
  $K_X^n>0$. Then $(X,L)$ is slope stable with respect to smooth subschemes
if
  \begin{enumerate}
  \item $2\mu(X,L)L + nK_X$ is nef, or
  \item $-2\mu(X,L)L - nK_X$ is nef.
\end{enumerate}
Moreover, for any divisor $G$ there is a $\delta_0>0$ such that if $0\le
\delta<\delta_0$ and $L=K_X(\delta G)$ is ample then $(X,L)$ is slope stable with respect to smooth subschemes.
\end{thm}

\begin{rmk}\label{rmk:weinkove}
  Suppose $K_X$ is ample. Then there exists an open set around
  $-[c_1(X)]$ of classes which admit cscK metrics \cite{LeBS}, so one
  has slope-semistability for these classes.
  
  It is shown in \cite{We} that if
  $$-2\mu(X,L)L - (n-1)K_X$$
  is ample then the Mabuchi functional
  associated to the class $c_1(L)$ is bounded from below, confirming the
  second result of Theorem \ref{canonicalbundleimpliesstability}.
\end{rmk}

\begin{proof}
  Fix a $Z$ and suppose $0<x<c\le \epsilon(Z)$.  Let $$f(x)=
  2n!(n-1)![a_0a_1(x) - a_1a_0(x)].$$
  We will show that $\int_0^c f(x)
  dx\le0$ for all smooth subschemes $Z$ and all $0<c\le \epsilon(Z)$,
  which implies $\mu_c(\I_Z)<\mu(X,L)$ since $a_0'(x)<0$.  For the third
part we will show this holds as long as $\delta<\delta_0$ where
  $\delta_0$ will be chosen independently of $Z$ and $c$.
  
  For $x\in(0,\epsilon(Z))$, $L_x:=L(-xE)$ is nef, so as $K_{\widehat
    X}-K_X=(p-1)E$ is effective,
  \begin{eqnarray*}
    f(x) &=& -(L^n)L_x^{n-1}.K_{\widehat X}+ (K_X.L^{n-1})L_x^n \\
 &\le& -(L^n)L_x^{n-1}.K_X+ (K_X.L^{n-1})L_x^n \\
&=& L_x^{n-1}.(B-x(K_X.L^{n-1})E),
  \end{eqnarray*}
  where
\begin{eqnarray*}
  B &:=& (K_X.L^{n-1})L - (L^n)K_X \\
&=& (K_X.L^{n-1})(L-K_X) -  ((L-K_X).L^{n-1})K_X\\
&=& \delta (K_X.L^{n-1}) G - \delta (G.L^{n-1})K_X.
\end{eqnarray*}
Notice that $B.L^{n-1}=0$.  Now, if $-B =
\frac{L^n}{n}(2\mu(X,L)L+nK_X)$ is nef then, as $L_x$ is nef,  $f(x)\le 0$, which proves (1).  When $n=1$ (so $X$ is a smooth curve), $B$ is
numerically trivial so $f(x)\le 0$ and we are done.  So we suppose
that $n\ge 2$.

As $B=O(\delta)$ we certainly have $f(x)\le 0$ for $\delta$
sufficiently small for any fixed value of $x$. However, since such a
choice of $\delta$ is not uniform in $x$, we integrate:
\begin{equation}
  \label{integraloff}
  \int_0^c f(x)dx = I_1-(K_X.L^{n-1})I_2,
\end{equation}
where $I_1 = \int_0^c L_x^{n-1}.B dx$ and $I_2=\int_0^c xL_x^{n-1}.E
dx$.  Then
\begin{eqnarray*}
 L_x^{n-1}.B &=& L^{n-1}.B+(L_x-L).\sum_{j=0}^{n-2} L^j.L_x^{n-2-j}.B\\
&=& -x \sum_{j=0}^{n-2} L^j.L_x^{n-2-j}.E.B, \qquad\text{as $L^{n-1}.B=0$.}
 \end{eqnarray*}
 We claim that for any $a$ and $b$,
\begin{equation*}
\sum_{j=0}^{n-2} \int_0^c xa^j (a-xb)^{n-2-j} dx = \frac{c^2}{n} \sum_{j=0}^{n-2}
(j+1) a^{j}(a-cb)^{n-2-j},
\end{equation*}
which can be shown by comparing the coefficient of $c^{n-j}$ on both
sides for $j=0,\ldots,n-2$ and using the identity
$$
\sum_{j=0}^i (j+1) \binom{n-2-j}{i-j} = \frac{n}{n-i} \sum_{j=0}^i
\binom{n-2-j}{i-j}\quad{\text{for $i=0,\ldots,n-2$}}.$$
Hence
\begin{equation}\label{I1}
 I_1 = \int_0^c L_x^{n-1}.B dx \le 
 c^2 \sum_{j=0}^{n-2} L^j.L_c^{n-2-j}.E.B.
\end{equation}
Similarly as
\begin{eqnarray*}
  \int_0^c x (a-xb)^{n-1}dx = \frac{c^2}{n(n+1)}\sum_{j=0}^{n-1}(n-j)a^j (a-cb)^{n-1-j},
\end{eqnarray*}
\begin{equation}\label{I2}
I_2=\int_0^c xL_x^{n-1}.E dx \ge \frac{c^2}{n(n+1)} \sum_{j=0}^{n-2} L^{j+1}.L_c^{n-2-j}.E.
\end{equation}
Putting (\ref{integraloff}), (\ref{I1}), (\ref{I2}) together
\begin{eqnarray*}
  \int_0^c f(x) dx &\le& I_1- (K_X.L^{n-1})I_2 \\
&\le&  -c^2\left(B + \frac{(K_X.L^{n-1})}{n(n+1)}L\right).\sum_{j=0}^{n-2}L^j.L_c^{n-2-j}.E.\\
\end{eqnarray*}

Recall that $L$ and $L_c$ are nef classes.  So it is now sufficient to prove
that $B+\frac{(K_X.L^{n-1})}{n(n+1)}L$ is also nef. But
$$
B+\frac{(K_X.L^{n-1})}{n(n+1)}L = \frac{L^n}{n}\left(-2\mu(X,L)L- n
  K_X - \frac{2\mu(X,L)}{n(n+1)}L\right).$$
As $\mu(X,L)\le 0$, this
is nef when $-2\mu(X,L)L - nK_X$ is, proving (2).

To prove the third part we must show that $\int_0^c f(x)\le 0$ uniformly with
respect to $\delta$.  Notice that the statement of the theorem is
unchanged if we scale $G$ by some positive number.  So without loss of
generality we suppose that $K_X(G)$ is ample.  Now
\begin{multline*}
    \quad B+\frac{(K_X.L^{n-1})}{n(n+1)}L\,=\,\frac{(K_X.L^{n-1})}{2n(n+1)}
    (L + (2n(n+1)+1)\delta G)  \\
    +\left(\frac{(K_X.L^{n-1})}{2n(n+1)}-\delta(G.L^{n-1})\right)K_X.\quad
\end{multline*}
For positive $\delta$ sufficiently small the line bundle
$$L+(2n(n+1) +1)\delta G=K_X + (2n(n+1)+2)\delta G$$
is ample, for it lies
on the line between $K_X$ (which is nef) and $K_X(G)$ (which is ample).
Moreover
$$
\frac{(K_X.L^{n-1})}{2n(n+1)}-\delta(G.L^{n-1}) =
\frac{K_X^n}{2n(n+1)} + O(\delta)$$
is positive for $\delta$
sufficiently small, since $K_X^n>0$.  Hence
$B+\frac{(K_X.L^{n-1})}{n(n+1)}L$ is nef for $\delta$ sufficiently
small, and the proof is complete.
\end{proof}

\subsection{Slope stability of smooth curves}\label{sec:slopestabilitysmoothcurves} \vskip 5pt
\noindent Since smooth curves always have cscK metrics, they should be stable.  We give a direct proof that they are slope (poly)stable:

\begin{thm}\label{smoothcurvesareKstable} Any smooth polarised curve $(\Sigma,L)$ of genus $g$ is  slope stable if $g\ge1$ and strictly slope polystable if
$g=0$.
\end{thm}

\begin{proof}
  Any nonempty subscheme $Z$ is a divisor of degree $d>0$, so
  $$
  \chi(L^{\otimes k}\otimes\I_Z^{xk})=k\deg L-xdk+1-g
  $$
  which shows that $\tilde a_0(x)=xd$ and $\tilde a_1(x)=0$. Thus
  $\mu_c(\O_Z)= \frac{cd}{c^2d}=\frac1{c}>0\ge\frac{1-g}{\deg
    L}=\mu(X)$ for $g\ge1$, proving slope stability.
    
    For $g=0$, $c$
  may take values up to and \emph{including} $\epsilon(Z)=\deg L/d$, since
  $L^{\otimes d}\otimes\I_Z^{\deg L}=\O_{\PP^1}(d\deg L-d\deg L)=\O_{\PP^1}$ is
globally generated.
Thus $\mu_c(\O_Z)\ge\frac{d}{\deg L}\ge\frac1{\deg L}=\mu(X)$ with equality
(strict semistability) only for $d=1$, \emph{i.e.}\ $Z$ a single point, and
$c=\epsilon(Z)$. Since the deformation to
  the normal cone of a single point on $\PP^1$ blows down to $\PP^1\times\C$
  (with a nontrivial $\C^\times$-action) from which the relevant line bundle
  $\L_\epsilon$ pulls back, we find $\PP^1$ is in fact slope polystable.
  \end{proof} 

\begin{rmk}
  In Corollary 6.7 of \cite{RT} it is shown that, for smooth curves,
  slope (semi/poly)stability is equivalent to K-(semi/poly)\-stability.
  Thus smooth curves are algebraically K-stable for $g\ge 1$ and
  algebraically K-polystable if $g=0$.
\end{rmk}

\subsection{Projective bundles}\label{proj}\ \vskip 5pt
Fix a polarised manifold $(B,\O_B(1))$ of dimension $b$, and
let $E$ be a vector bundle on $B$ with $r+1:=\rk E\ge 2$.  We show
that the stability of $\PP(E)$ is related to slope stability of the
bundle $E$ (as defined in \ref{sec:slopestabbundle}) and slope
stability of the base $B$.  Let $n=\dim \PP(E)=b+r$ and
$$
L_m= \O_{\PP(E)}(1) \otimes \O_{B}(m),$$
which is ample for $m$
sufficiently large.

\begin{thm}\label{thm:slopestableprojectivebundles} 
  If $(\PP(E),L_m)$ is slope semistable for all $m\gg 0$ then $E$ is
  a slope semistable vector bundle and $(B,\O_B(1))$ is a slope semistable
  manifold. Moreover, there is an $m_0$ which depends only on $E$ and $(B,\O_B(1))$
  such that if $(\PP(E),L_m)$ is slope semistable for some $m\ge m_0$ then
  $E$ is a slope semistable vector bundle.
  
  Thus if $E$ is a strictly slope unstable bundle or if $(B,\O_B(1))$
  is a strictly slope unstable manifold, then $\PP(E)$ does not admit
  a cscK metric in $[c_1(L_m)]$ for $m\gg 0$.
\end{thm}

For bundles (of any rank) over curves, we get stronger results.

\begin{thm}\label{thm:slopestableprojectivebundlesovercurves}
  Suppose $B$ is a smooth curve of genus $g\ge 1$ and that $L_m$ is ample.
  If $(\PP(E),L_m)$ is slope (semi/poly)stable then $E$ is slope
  (semi/poly)stable.
  
  If $E$ is polystable then $\PP(E)$ has a cscK metric in every
  K\"ahler class. Conversely if $E$ is strictly unstable then $\PP(E)$
  does not admit a cscK metric in any rational K\"ahler class.
  Finally, if $E$ is not polystable then $\PP(E)$ is not algebraically
  K-polystable.
\end{thm}

The proofs appear after a calculation of the relevant slopes and
Seshadri constants.  It is well known that if $E$ is polystable and
$B$ is a curve then $\PP(E)$ admits a cscK metric in every K\"ahler
class \cite{BdB}.  So Theorem
\ref{thm:slopestableprojectivebundlesovercurves} gives an almost
complete converse.  If, as expected, a cscK metric implies algebraic
K-polystability then Theorem
\ref{thm:slopestableprojectivebundlesovercurves} would be a full
converse.  Moreover it would imply that slope polystability is
equivalent to algebraic K-polystability for projective bundles over
curves of genus $g\ge 1$.

There is also a partial converse to Theorem
\ref{thm:slopestableprojectivebundles}.  Suppose that $E$ is slope
stable, and $B$ is a manifold with $\aut(B)=0$ and a cscK metric in
$c_1(\O_B(1))$.  Then there exists a cscK metric on $\PP(E)$ in
$c_1(L_m)$ for $m$ sufficiently large \cite{Ho}.

In the rank$\,E=2$, dim\,$B=1$ case,
it is known that if a ruled surface $\PP(E)$ has a cscK metric in any
class then $E$ is a polystable bundle.  This is proved by \cite{BdB}
in the scalar-flat case, by \cite{LeB} in the case that $g\ge 2$,
$-\int_X c_1(K_X).c_1(L)<0$, and \cite{ATF} in general.

The stability of ruled surfaces has also been studied by Morrison \cite{Mo}.
If $E$ is unstable
then $\PP(E)$ is Chow unstable with respect to what he calls ``good''
polarisations (in particular $(\PP(E), L_m^{\otimes k})$ is Chow unstable for
$k\gg 0$).  By \cite{Do2} this implies that $\PP(E)$ does not have a cscK
metric in any class.  Morrison also shows that if $E$ is stable then
for suitable $m$,  $(\PP(E),L_m)$ is Chow stable, and he
conjectures that this holds for $(\PP(E),L_m^{\otimes k})$ with $k\gg
0$.  Since there exists a cscK metric in $c_1(L_m)$, this
conjecture follows from \cite{Do2} when $g\ge 2$ and $E$ is simple.

\begin{rmk}\label{ext1}
  Suppose that $E\to B$ has a subbundle $F$. Let $\X$ be the degeneration
  to the normal cone of $\PP(F)\subset(\PP(E),L_m)$ with
  $c=\epsilon(\PP(F))=1$.  Then $\L_c$ is only semi-ample but not ample,
  and contracts a component of the central fibre $\X_0$ (test
  configurations with semi-ample polarisation are studied in
  Proposition 5.1 of \cite{RT}).  This contraction is a test configuration which is the projectivisation
  of the degeneration of bundles taking the extension
$$
0\to F\to E\to G\to0, \quad\text{defined by }\ e\in\mathrm{Ext}^1(G,F),
$$
to the direct sum $F\oplus G$ (via the family of extensions $\lambda
e,\ \lambda\in\C$). If $e=0$ (\emph{i.e.}\ $E=F\oplus G$ to begin with) then
we get a product degeneration.  We show below that if $F$ and $E$ have the same slope then the Futaki invariant is 0 (on curves, and to the
two top orders in $m$ for general $B$).  So we recover the usual notion of
polystability for bundles.
\end{rmk}

\subsubsection{Slope stability of vector bundles}\label{sec:slopestabbundle}

For brevity write
$\mu(B)=\mu(B,\O_B(1))$.  For any coherent sheaf $E$ on $B$ it is convenient
to define
$$\mu\_E=\frac{\deg E}{a_0^B(b-1)!\rk E}+ \mu(B),$$
where
$\chi(\O_B(k)) = a_0^B k^b + a_1^B k^{b-1} + O(k^{b-2})$.  Note that
this differs from the usual definition of slope for a sheaf. However for
any coherent subsheaf $F$,
$$\mu\_E-\mu\_F = \frac{1}{a_0^B(b-1)!}\left(\frac{\deg E}{\rk
    E}-\frac{\deg F}{\rk F}\right).$$
Thus $E$ is a slope stable
(resp.\ semistable) vector bundle if and only if $\mu\_F<\mu\_E$ (resp.\ 
$\mu\_F\le \mu\_E$) for all coherent subsheaves $F<E$.  And $E$ is
polystable if and only if it is a direct sum $E=\oplus F_i$ of slope
stable sheaves, with $\mu\_{F_i}=\mu\_E$ for all $i$.

\begin{lem}\label{lem:somefactsaboutslope}
  Let $E$ and $F$ be torsion free coherent sheaves on $B$. Then
  \begin{enumerate}
  \item $\chi(E\otimes\mathcal O_B(m)) = a_0^B \rk E (m^b + \mu\_E
    m^{b-1}) + O(m^{b-2})$, where the $O(m^{b-2})$ is understood to be
    zero when $b=\dim B=1$,
  \item $\mu\_{S^k E^*} = (1+k)\mu(B)-k\mu\_E$,
  \item $\mu\_{E\otimes F} = \mu\_E + \mu\_F - \mu(B)$,
  \item if $F<E$ and $E/F$ is also torsion free then $$(\rk E) \mu\_E = (\rk
    F) \mu\_F + (\rk(E/F))\mu\_{E/F}.$$
  \end{enumerate} 
\end{lem}

\begin{proof}
  From the definition of $\mu\_E$ and the Riemann-Roch theorem,
\begin{eqnarray*}
\chi(E\otimes \mathcal O_B(m)) &=& \int_B \ch(\mathcal O_B(m))\text{ch}(E)\text{Td}_B \\
&=& \int_B e^{mc_1(\mathcal O_B(1))}(\rk E+c_1(E)+ \cdots)\text{Td}_B\\
&=& a_0^B \rk E (m^b + \mu\_E  m^{b-1}) + O(m^{b-2}).
\end{eqnarray*}
Now as $E$ is torsion free, we can
calculate the degrees of $E$ and $S^k E^*$ by restricting to the set
where $E$ is locally free, since its complement has codimension
$\ge2$. We compute $\mu\_{S^k E^*}$ to be
\begin{eqnarray*}
 \frac{\deg S^k E^*}{a_0^B (b-1)! \rk S^k E} + \mu(B)
&=&-k\frac{\deg E}{a_0^B (b-1)! \rk E} + \mu(B)\\
&=&(1+k)\mu(B)-k\mu\_E,
\end{eqnarray*}
where the second equality follows from the splitting principle.  Also,
\begin{equation*}
    \mu\_{E\otimes F} = \frac{\rk F\deg E+\rk E\deg F}{a_0^B (b-1)! \rk E \rk F} + \mu(B) = \mu\_E + \mu\_F -\mu(B).
\end{equation*}
Finally if $F<E$ then comparing the $m^{b-1}$ terms in $\chi(E\otimes \O_B(m)) = \chi(F\otimes\O_B(m)) + \chi((E/F)\otimes \O_B(m))$ gives $(\rk E) \mu\_E = (\rk F) \mu\_F + (\rk(E/F))\mu\_{E/F}$.
\end{proof}

\subsubsection{Seshadri constants of projective subbundles}

For the rest of this section let $Z=\mathbb P(F)$, so the Seshadri
constant $\epsilon(\mathbb P(F),L_m)$ is defined as in
(\ref{def:seshadri}).

\begin{lem} \label{Groth}
For $E$ be a vector bundle over a curve $B$,
$\deg_B\!E^*=\deg_{\,\PP(E)}\!\O_{\PP(E)}(1).$
\end{lem}
\begin{proof}
Let $\omega$ denote $c_1(\O_{\PP(E)}(1))$ on $\PP(E)$.
The general Grothendieck formula $\sum_{i=0}^{r+1}\omega^{r+1-i}c_i(E)=0$ reduces over a curve to $-c_1(E)\omega^r=\omega^{r+1}$, whose left hand side
is $-\deg_B E$.
\end{proof}

\begin{prop}\label{prop:seshadri}
  
  There is an $m_0$ \emph{(}depending only on $E$ and the pair $(B,\O_B(1))$ such
  that for any $m\ge m_0$ and any \emph{saturated} subsheaf $F$ of the
  bundle $E$ \emph{(i.e.}\ $E/F$ is torsion free\emph{)} with $\mu\_F\ge
  \mu\_E$ we have $\epsilon=\epsilon(\PP(F),L_m)=1$.
  
  Suppose that $B$ is a curve, $F<E$ is saturated, $E/F$ is semistable
  and $\mu\_F\ge \mu\_E$.  Then for any $m$ such that $L_m$ is ample,
  $\epsilon(\PP(F),L_m)=1$ and the global sections of
  $L_m^{\otimes k}\otimes\I_{\PP(F)}^k$ generate $\I_{\PP(F)}^k$ for $k\gg0$.
\end{prop}

\begin{proof}
  Since $\mathbb P(F_p)\subset \mathbb P(E_p)$ is a linear subspace
  for any $p\in B$, and $L_m|_{\mathbb P(E_p)}=\O_{\PP(E_p)}(1)$,
  it follows that $\epsilon\le 1$.  To show $\epsilon$ is at least 1 it is
  sufficient to show that $L_m\otimes\I_{\mathbb P(F)}$ is generated by global
  sections.
  
  Let $G=E/F$.  As the set of quotients $G$ of $E$ with $\mu\_G\le\mu\_E$
  is bounded (\cite{HL} Lemma 1.7.9) there is an $m_0$ (depending only
  on $E$ and $(B,\O_B(1))$) such that for all $m\ge m_0$, $G^*(m)$ is
  globally generated and has no higher cohomology and $L_m$ is ample.

  
  Working on $\mathbb P(E)$, $\mathcal O_{\mathbb P(E)}(-1)$ is a
  subbundle of (the pullback of) $E$ giving a canonical element
  $u\in\Hom(\mathcal O_{\mathbb P(E)}(-1), G)$ obtained by composition
  with the projection from $E$ to $G$.  Thinking of $u$ as a section
  of $G\otimes \O_{\PP(E)}(1)$, its zero set is precisely $\mathbb
  P(F)$.
  
  Now turn to $\mathbb P=\mathbb P(H^0(L_m)^*)=\mathbb
  P(H^0(E^*(m))^*)$ and let $m\ge m_0$.  The exact sequence $0\to
  H^0(G^*(m))\to H^0(E^*(m))\to H^0(F^*(m))\to 0$ yields a canonical
  section $v$ of $H^0(G^*(m))^*\otimes\mathcal O_{\mathbb P}(1)$ whose
  zero set is $\mathbb P(H^0(F^*(m))^*)$.
  
  Since $G^*(m)$ is globally generated, $G(-m)$ injects into
  $H^0(G^*(m))^*$. Tensoring with $L_m$ shows
  that $G\otimes \O_{\PP(E)}(1)$ injects into $H^0(G^*(m))^*\otimes L_m$,
  and $u$ maps to $v$.  Thus $\mathbb P(F)$ is the intersection of
  $\mathbb P(E)$ with the subspace $\mathbb P(H^0(F^*(m))^*)$ of
  $\mathbb P$.  Hence $L_m\otimes \I_{\mathbb P(F)}$ is generated by
  global sections, so $\epsilon= 1$ as claimed.
  
  Now suppose that $B$ is a curve, $m$ is chosen so that $L_m$ is ample,
  $\mu\_F\ge \mu\_E$, $F$ is destabilising and saturated, and
  $G=E/F$ is semistable. Since $F$ is saturated, $G$ is torsion free, so
  both are locally free since $B$ is a curve.  Then $\mu\_E\ge \mu\_G$, so
  $$
  \deg(G^*\otimes \O_B(m))\ge \deg(E^*\otimes \O_B(m))=\deg L_m>0,
  $$
  by Lemma \ref{Groth}.
  Thus $G^*\otimes \O_B(m)$ is a semistable bundle of positive degree on
  a curve $B$, so it is ample, \emph{i.e.}\ 
  $\O_{\PP(G)}(1)\otimes \O_B(m)$ is ample (\cite{La} 6.4.11).  So
  for $k\gg 0$, $\O_{\PP(G)}(k)\otimes \O_B(km)$ is globally generated
  and thus so is its pushdown $S^k G^*\otimes \O_B(km)$.
  
  Now $\pi_*(L^{\otimes k}_m\otimes\I_{\PP(F)}^k)= S^kG^*\otimes \O_B(mk)\ \big(<
  S^k E^*\otimes\O_B(mk)\big)$ generates $\I_{\PP(F)}^k$ on each fibre, so
  the global sections of $L_m^{\otimes k}\otimes\I_{\PP(F)}^k$ generate $\I_{\PP(F)}^k$
  as claimed.  From the definition of the Seshadri constant we again get
  that $\epsilon= 1$.
\end{proof}

\subsubsection{Slope of projective bundles}\

In calculating the quantities $\mu(\PP(E),L_m)$ and $\mu(\O_{\PP(F)},L_m)$ it is convenient to make the change of variables 
\begin{equation}
  \widetilde{m} = m + \frac{1}{b}(\mu(B)-\mu\_E).\label{eq:definitionoftildem}
\end{equation}
 (The reader may prefer to assume that $\deg E=0$, in which case $\widetilde{m}=m$.)   We write $\mu(\PP^r):=\mu(\PP^r,\O_{\PP^r}(1)) = r(r+1)/2$. 

\begin{lem}\label{lem:eulerprojbundle}
  Let $\chi(\mathbb P(E),L_m^{\otimes k}) = a_0 k^n + a_1 k^{n-1} + O(k^{n-2})$.
  Then $a_0$ and $a_1$ are polynomials in $m$.  In fact if $\widetilde{m}$
  is defined as in \eqref{eq:definitionoftildem} then
\begin{eqnarray*}
  a_0 &=&  \frac{a_0^B}{r!}\widetilde{m}^b +O(\widetilde{m}^{b-2}), \qquad\text{and}\\
  a_1 &=&  \frac{a_0^B}{r!}\left(\mu(\PP^r)\widetilde{m}^b + \mu(B)\widetilde{m}^{b-1}\right)
 +O(\widetilde{m}^{b-2}),
\end{eqnarray*}
where if $\dim B=1$ we interpret $O(\widetilde{m}^{b-2})$ as being zero.
Moreover the $O(\widetilde{m}^{b-2})$ terms depend only on $(B,\O_B(1))$ and the
Chern classes of $E$.
\end{lem} 
\begin{proof}
Let $\pi\colon\PP(E)\to B$ be the
  projection.   As $L_m$ is relatively ample, for $k\gg 0$,
\begin{eqnarray*}
  \chi(L_m^{\otimes k})\!\!&=&\!\!\chi(\pi_*(L_m^{\otimes k}))\\
&=&\!\!\chi(S^k E^* \otimes \mathcal O_B(mk))\\
&=&\!\! a_0^B \rk S^k E \cdot \left(m^b k^b + \mu\_{S^k E^*}
  m^{b-1}k^{b-1}\right) +O(m^{b-2}) \\
&=&\!\! a_0^B \rk S^k E \cdot \left( m^b k^b + [(1+k)\mu(B) -
  k\mu\_E]m^{b-1}k^{b-1}\right) \\
&&+O(m^{b-2}),
\end{eqnarray*}
where in the last line we have used (\ref{lem:somefactsaboutslope} (2)) and the $O(m^{b-2})$ term is zero if $b=1$. Now the rank term is
$$\rk S^k E = \binom{r+k}{k} = \frac{1}{r!}\left[ k^r + \mu(\mathbb
  P^r)k^{r-1} + O(k^{r-2})\right].$$
Expanding and taking the $k^n$
and $k^{n-1}$ terms gives
\begin{eqnarray*}
  a_0 &=&  \frac{a_0^B}{r!}\left (m^b + [\mu(B)-\mu\_E]m^{b-1} \right)+ O(m^{b-2}),\\
  a_1 &=&  \frac{a_0^B}{r!}\left( \mu(\mathbb P^r)m^b +
  \mu(\mathbb P^r)[\mu(B)-\mu\_E]m^{b-1} + \mu(B)m^{b-1}\right) \\
&&+ O(m^{b-2}),
\end{eqnarray*}
and the change of variables from $m$ to $\widetilde{m}$ gives the
expressions in the statement of the lemma.  As the Chern character of
$S^kE$ depends only on $k$ and the Chern classes of $E$ (\cite{Ha}
Appendix A3) we see that the $O(m^{b-2})$ terms (and hence the  $O(\widetilde{m}^{b-2})$ terms) depend only on
$(B,\O_B(1))$ and the Chern classes of $E$.
\end{proof}

\begin{lem}\label{lem:wtsofprojsubbundle}
  Let $F$ be a saturated coherent subsheaf of $E$.  Define
  $\alpha_i(x)$ for $\PP(F)\subset\PP(E)$ as in Proposition
  \ref{quotientslopeusingalpha}, and let $\widetilde{m}$ be defined as
  in \eqref{eq:definitionoftildem}. Then,
\begin{eqnarray*}
  \int_0^1 (1-x) \alpha_1(x) dx &= \frac{a_0^B(s+1)}{(r+1)!}
  \left[\widetilde{m}^b + \frac{1}{r+2}(\mu\_E-\mu\_F)\widetilde{m}^{b-1}\right]
  + O(\widetilde{m}^{b-2}),
\end{eqnarray*}
\begin{multline*} 
\!\!\!\! \int_0^1 \left((1-x)\alpha_2(x) + \frac{\alpha_1(0)}{2}\right) dx = \\
\frac{a_0^B(s+1)}{2(r+1)!}\left(2\mu(\PP^r) \widetilde{m}^b + \left[2\mu(B) +
(r+1)(\mu\_E-\mu\_F)\right]\widetilde{m}^{b-1}\right) + O(\widetilde{m}^{b-2}),
\end{multline*}
where $O(\widetilde{m}^{b-2})$ is understood to be zero if $\dim B=1$.  Both
expressions are polynomials in $\widetilde{m}$, and the $O(\widetilde{m}^{b-2})$
terms depend only on $(B,\O_B(1))$ and the Chern classes of $E$ and
$F$.
\end{lem}
\begin{proof}
  Let $F$ have rank $s+1$ and $G=E/F$ have rank $t+1$, so 
  \begin{equation}
    \label{eq:rst}
    s+t+2=r+1.
  \end{equation}
  Since $F$ is saturated, $G$ is torsion free, as is $F$ since it sits
  inside a locally free sheaf $E$. Thus $E$, $F$ and $G$ are locally
  free on an open set $U$ whose complement has codimension at least
  $2$, on which their first Chern classes (and those of their
  symmetric powers) can be calculated. Since $G$ is locally free on
  $U$, $F\into E$ has constant rank so $\mathbb P(F|_U)$ sits inside
  $\mathbb P(E|_U)$ as a smooth submanifold with normal bundle
$$
\nu = \nu_{\mathbb P(F|_U)}=\pi^* G \otimes \mathcal O_{\mathbb
  P(F)}(1),
  $$
where, by abuse of notation, we let $\pi=\pi|_{\PP(F)}$.
As the complement of $V=\mathbb P(E|_U)$ also has
codimension at least $2$, we can calculate $\alpha_1(x)$ and
$\alpha_2(x)$ on $V$. Then for $0<x<1$,
$$
\pi_*(L^{\otimes k}_m\otimes S^{xk} \nu^*) = S^{xk} G^* \otimes S^{(1-x)k
}F^* \otimes \mathcal O_B(mk).$$
Furthermore for
$0<x\ll 1$ the higher cohomology of both
$L_m^{\otimes k}\otimes \I_{\PP(F)}^{xk}$ and $L_m^{\otimes k}\otimes \I_{\PP(F)}^{xk+1}$
vanish for $k\gg 0$ and hence the same is true for $L_m^{\otimes k}\otimes
\I_{\PP(F)}^{xk}/\I_{\PP(F)}^{xk+1}$, so
\begin{equation} \label{forms}
\chi(S^{xk} G^* \otimes S^{(1-x)k}F^* \otimes \mathcal O_B(mk)) = \alpha_1(x) k^{n-1} + \alpha_2(x)k^{n-2} + O(k^{n-3}).
\end{equation}
Now let $R=\rk S^{(1-x)k} F \cdot
\rk S^{xk} G$, which equals
\begin{align*}
\frac{1}{s!t!}&\left[(1-x)^s k^s +
  \frac{s(s+1)}{2}(1-x)^{s-1}k^{s-1}+\cdots \right]\cdot\\
&\left[x^tk^t +
  \frac{t(t+1)}{2} x^{t-1}k^{t-1}+\cdots \right]\\
&=\frac{1}{s!t!}\left((1-x)^sx^tk^{r-1} + \delta(x) k^{r-2} + O(k^{r-3})\right),
\end{align*}
where $2\delta(x) = s(s+1)(1-x)^{s-1}x^t + t(t+1)(1-x)^{s}x^{t-1}$.
Notice that this holds even if $s$ or $t$ are zero, for then the
$k^{s-1}$ or $k^{t-1}$ terms vanish.  Much calculation with Lemma
\ref{lem:somefactsaboutslope} computes
\begin{eqnarray*}
&& \hspace{-1cm} \chi(S^{xk} G^* \otimes S^{(1-x)k}F^*\otimes \mathcal O_B(mk))\\
&=&a_0^B R\big(m^bk^b +  [(1+k)\mu(B) - kx\mu\_G-k(1-x)\mu\_F]m^{b-1}k^{b-1}\big)\\
&&+O(m^{b-2}) \\
&=&a_0^B R((m^b + \mu(B)m^{b-1})k^b +  \mu(B)m^{b-1}k^{b-1}) \\
&&- a_0^B R(x\mu\_G+(1-x)\mu\_F)m^{b-1}k^{b} + O(m^{b-2}).
\end{eqnarray*}
Now $m^b + \mu(B)m^{b-1} = \widetilde{m}^b+\mu\_E\widetilde{m}^{b-1} + O(\widetilde{m}^{b-2})$, and
\begin{eqnarray*}
\gamma(x)&:=&  \mu\_E -(x\mu\_G + (1-x)\mu\_F) \\ 
&=& \mu\_E - \frac{x}{t+1}((r+1)\mu\_E - (s+1)\mu\_F) - (1-x)\mu\_F \\
&=& (\mu\_E-\mu\_F)\left(1-\frac{x(r+1)}{t+1}\right),
\end{eqnarray*}
where the last line uses (\ref{eq:rst}).  Thus $\chi(S^{xk} G^* \otimes
S^{(1-x)k }F^* \otimes \mathcal O_B(mk))$ is
$$a_0^B R (\widetilde{m}^bk^b + \gamma(x)\widetilde{m}^{b-1}k^b + \mu(B)\widetilde{m}^{b-1}k^{b-1})+O(\widetilde{m}^{b-2}).$$
Now $\alpha_1(x)$ and
$\alpha_2(x)$ (\ref{forms}) are polynomials in $x$ so extend uniquely
from $0<x\ll 1 $ to all of $\mathbb R$, and the above shows that
\begin{eqnarray*}
  \alpha_1(x) &=& \frac{a_0^B}{s!t!}(1-x)^s x^t(\widetilde{m}^b + \gamma(x)\widetilde{m}^{b-1})
+ O(\widetilde{m}^{b-2}),
\quad\text{and}\\
  \alpha_2(x) &=& \frac{a_0^B}{s!t!}\delta(x)(\widetilde{m}^b + \gamma(x)\widetilde{m}^{b-1})
 +\frac{a_0^B}{s!t!}\mu(B)(1-x)^s x^t\widetilde{m}^{b-1} \\
&&+ O(\widetilde{m}^{b-2}).
\end{eqnarray*}

To calculate the required integrals of the $\alpha_i(x)$ one has to
consider four cases, depending on whether $s$ or $t$ vanish.  In all
four cases, repeated applications of the identity $\int_0^1 (1-x)^s
x^tdx = \frac{s!t!}{(s+t+1)!}$ give the formula in the statement of the
Lemma.  These expressions depend only on $(B,\O_B(1))$ and the Chern
characters of the symmetric powers of $E$ and $F$, and thus only on
$(B,\O_B(1))$ and the Chern classes of $E$ and $F$ by (\cite{Ha} Appendix
A3).
\end{proof}


\begin{prop}\label{prop:slopeprojectivebundles}\
  Let $F$ be a saturated coherent subsheaf of $E$ and suppose that
  $\epsilon(\PP(F),L_m)=1$.   Then $\mu_{1}(\O_{\PP(F)},L_m)-
  \mu(\PP(E),L_m)$ equals
\begin{equation*}
C\left((\mu\_E-\mu\_F)\big[(r+1)\widetilde{m}^{2b-1} -\mu(B)\widetilde{m}^{2b-2}\big]
+ O(\widetilde{m}^{2b-3})\right), \end{equation*}
where $C=C(\widetilde{m})$ is positive.  Here the $O(\widetilde{m}^{2b-3})$ term is understood to be zero if $B$ is a curve.  Moreover the $O(\widetilde{m}^{2b-3})$ terms
depend only on $(B,\O_B(1))$ and the Chern classes of $F$ and $E$.
\end{prop}

\begin{proof}
Using the expressions in Lemmas \ref{lem:eulerprojbundle} and 
 \ref{lem:wtsofprojsubbundle} gives
 \begin{multline*}
   a_0\int_0^1 \left((1-x)\alpha_2(x) + \frac{\alpha_1(0)}{2}\right) dx - a_1\int_0^1 (1-x)\alpha_1(x)dx \\
= \frac{(a_0^B)^2(s+1)}{(r+2)!r!}(\mu\_E-\mu\_F)[(r+1)\widetilde{m}^{2b-1} -\mu(B)\widetilde{m}^{2b-2}]+ O(\widetilde{m}^{2b-3}).
 \end{multline*}
Thus $ \mu_1(\O_{\PP(F)},L_m) - \mu(\PP(E),L_m)$ equals
\begin{eqnarray*}
C(\widetilde{m})\left((\mu\_E-\mu\_F)[(r+1)\widetilde{m}^{2b-1} - \mu(B)\widetilde{m}^{2b-2}]+ O(\widetilde{m}^{2b-3})\right), \end{eqnarray*}
where 
$$ C(\widetilde{m}) = \frac{(a_0^B)^2(s+1)}{(r+2)!r!a_0\int_0^1(1-x)\alpha_1(x)dx}\,,$$
which is positive.
As the $O(\widetilde{m}^{b-2})$ terms of $a_0$ and $a_1$, as well as
$\alpha_1(x)$ and $\alpha_2(x)$, depend only on $(B,\O_B(1))$ and the
Chern classes of $F$ and $E$ so does the $O(\widetilde{m}^{2b-3})$ term above.
\end{proof}

\begin{proof}[Proof of Theorem \ref{thm:slopestableprojectivebundlesovercurves}]
  If $E$ is not slope stable (resp.\ strictly unstable) then there is
  a maximally destabilising subsheaf of $E^*$ which is saturated and
  so locally free. Call its dual $G$ and let $F$ be the kernel of
  $E\to G\to0$. Then $F<E$ is saturated and locally free, $G$
  is semistable, and $\mu\_F\ge \mu\_E$ (resp.\ $\mu\_F>\mu\_E$).
  Therefore by Proposition \ref{prop:seshadri} $\epsilon(\PP(F),L_m)=1$ and
  the global sections of $L_m^k\otimes\I_{\PP(F)}^k$ saturate $\I_{\PP(F)}^k$
  for $k\gg0$.
  
  As $\deg L_m=a_0^B(r+1)\widetilde{m}>0$ (Lemma \ref{Groth}) we
  have $\widetilde{m}>0$.  And as $g\ge 1$, $\mu(B)\le 0$, so
  $(r+1)\widetilde{m}-\mu(B)>0$.  From Proposition \ref{prop:slopeprojectivebundles},
  $$\mu_{1}(\O_{\PP(F)},L_m)- \mu(\PP(E),L_m) =
  C(\mu\_E-\mu\_F)[(r+1)\widetilde{m}- \mu(B)],$$
  where $C>0$.  Thus if $E$
  is not slope (semi)stable then $\PP(E)$ is not slope (semi)stable.

  Finally suppose that $E$ is not polystable.  Then there is a subbundle
  $F$ with either $\mu\_F>\mu\_E$, which we have already dealt with, or
  $\mu\_F=\mu\_E$ and \emph{$F$ is not a direct summand}. The
  degeneration to the normal cone of $\PP(F)$ with $c=\epsilon=1$
  gives a test configuration with zero Futaki invariant whose central
fibre is $\PP(F\oplus
  E/F)$ (\ref{ext1}). This cannot be a product configuration since the central
  fibre is not isomorphic to $\PP(E)$.  So $(\PP(E),L_m)$ is not slope
  polystable.
\end{proof} \vskip 5pt

We could similarly now prove the first part of Theorem
\ref{thm:slopestableprojectivebundles}, but to prove all of it we first
calculate the slope and Seshadri constant of $\PP(E|_{B'})$.

\begin{lem}
\label{lem:asymptoticrroch}
Let $C$ and $D$ be torsion free sheaves on $B$ and suppose that
$\mu\_C=\mu(B)$. Then
$$
\chi(C\otimes D\otimes \mathcal O_B(mk)) = \rk C \cdot
\chi(D\otimes\mathcal O_B(mk)) + O(m^{b-2}).$$
  \end{lem}
  
  \begin{proof}
    By Lemma \ref{lem:somefactsaboutslope} (3) the hypotheses imply that
    $\mu\_{C\otimes D} = \mu\_D$. Now apply Lemma \ref{lem:somefactsaboutslope}
    (1) twice.
  \end{proof}

  \begin{prop} \label{b'} Let $B'$ be a subscheme of $B$.  Then for
    $m\gg 0$, $\epsilon(\PP(E|_{B'}),L_m) \ge
    m\epsilon(B',\O_B(1))+O(m^0)$, and
  $$
  \mu_{cm}(\I_{\PP(E|_{B'})},L_m) -\mu(\PP(E),L_m)=
  \frac1m[\mu_{c}(\I_{B'},\O_B(1))-\mu(B)]+O(m^{-2}).
  $$
\end{prop}

\begin{proof}
Pick an integer $u$ so that $E^*(u)$ is globally generated; then so is
$S^k(E^*(u)) =S^k E^* \otimes \mathcal O_B(ku)$ for all $k$.  We first
show that 
$$\epsilon(\PP(E|_{B'}),L_m) \ge(m-u)\epsilon(B',\mathcal
O_B(1))=m\epsilon(B')+O(m^0).$$
 By the definition of the Seshadri
constant, if $c<(m-u)\epsilon(B',\mathcal O_B(1))$ then $
\mathcal O_B((m-u)k)\otimes\I_{B'}^{ck}$ is globally generated for $k$
sufficiently large.  Hence for $k\gg 0$ the sheaf
  \begin{eqnarray*}
    \pi_*(L_m^{\otimes k} \otimes\I_{\PP(E|_{B'})}^{ck}) &=& 
\mathcal O_B(mk) \otimes S^k E^*\otimes \I_{B'}^{ck}  \\
&=& \mathcal O_B(mk-uk) \otimes S^k(E^*(u))\otimes \I_{B'}^{ck}
  \end{eqnarray*}
  is also globally generated, and thus so is
  $L_m^{\otimes k} \otimes\I_{\PP(E|_{B'})}^{ck}$, because $L_m^{\otimes k}$ is globally generated
  along the fibres.  This implies that $c\le
  \epsilon(\PP(E|_{B'}),L_m)$, so $(m-u)\epsilon(B',\mathcal O_B(1))\le
  \epsilon(\PP(E|_{B'}),L_m)$.
  
  Now we calculate the slope of $\PP(E|_{B'})$.  Since we are
  interested in $m\gg 0$, we may twist $E$ by some power of $\mathcal
  O_B(1)$ to assume, without loss of generality, that
  $\mu\_E=\mu(B)$ (\emph{i.e.}\ $\deg E=0$).  This power may not be integral, but
  that does not affect the purely
  numerical argument below; we just have to allow rational $m$. Let
  \begin{eqnarray*}
\chi\_{\PP(E)}(L_m^{\otimes k} \otimes\I_{\PP(E|_{B'})}^{xk}) &=& a_0(x) k^n +a_1(x) k^{n-1} + O(k^{n-2}), \quad\text{and} \\
\chi\_B(\mathcal O_B(k)\otimes\I_{B'}^{xk}) &=& b_0(x) k^n +b_1(x) k^{n-1} +
O(k^{n-2}).
  \end{eqnarray*}
  Fix $x<\epsilon(\PP(E|_{B'}),L_m)$ and suppose $k\gg 0$.  Then
  $$\pi_*(L_m^{\otimes k} \otimes\I_{\PP(E|_{B'})}^{xk})=\I_{B'}^{xk} \otimes \mathcal
  O_B(mk) \otimes S^k E^*,$$
  with the higher
  pushdowns zero.  From Lemma \ref{lem:somefactsaboutslope} (2), we
  have that $\mu\_{S^k E}=\mu(B)$ for all $k$, so Lemma
  \ref{lem:asymptoticrroch} yields
  \begin{eqnarray*}
\chi(L_m^{\otimes k} \otimes\I_{\PP(E|_{B'})}^{xk}) &=& \chi(\I_{B'}^{xk} \otimes
  \mathcal O_B(mk) \otimes S^k E^*) \nonumber\\
&=& \rk S^k E \cdot \chi(\I_{B'}^{xk} \otimes \mathcal O_B(mk))+ O(m^{b-2})\nonumber\\
&=& \rk S^k E \cdot( b_0(x/m) m^bk^b + b_1(x/m)m^{b-1}k^{b-1}) \\
&&+ O(m^{b-2}).
  \nonumber
\end{eqnarray*}
Now
$$\rk S^k E = \binom{r+k}{r} = \frac{1}{r!}\left( k^r + \mu(\mathbb
  P^r) k^{r-1} + \cdots\right),$$
where $\mu(\mathbb P^r) =\mu(\mathbb
P^r,\mathcal O_{\mathbb P^r}(1))=r(r+1)/2$.  Thus
\begin{eqnarray*}
  a_0(x) &=& \frac{1}{r!}b_0(x/m)m^b + O(m^{b-2}), \qquad \text{and} \\ 
  a_1(x) &=& \frac{1}{r!}\left[\mu(\mathbb P^r)b_0(x/m)m^b+
  b_1(x/m) m^{b-1}\right] + O(m^{b-2}).
\end{eqnarray*}
Hence 
\begin{eqnarray*} 
\mu_{mc}(\I_{\PP(E|_{B'})},L_m)\ &=& \frac{\int_0^{mc} \left(a_1(x) +
  \frac{a_0(x)'}{2} \right) dx}{\int_0^{mc} a_0(x) dx}\\
&=& \frac{\int_0^c \left(a_1(mx) + \frac{a_0'(mx)}{2}\right) dx}{\int_0^c
  a_0(mx) dx} \\
&=& \frac{\int_0^{c} \left(\mu(\mathbb P^r)m^b b_0(x) + m^{b-1}(b_1(x) +
  \frac{b_0(x)'}2) \right)dx + O(m^{b-2})}{\int_0^cm^bb_0(x) dx + O(m^{b-2})} \\
&=& \mu(\mathbb P^r) +\frac{1}{m}\mu_c(\I_{B'},\mathcal
  O_B(1))+ O(m^{-2}).
\end{eqnarray*}
On the other hand, by assuming $\deg E=0$ we have $\widetilde{m} =m$
so from Lemma (\ref{lem:eulerprojbundle}),
$$\mu(\PP(E),L_m) = \frac{a_1}{a_0}=\mu(\PP^r) + \frac{1}{m}\mu(B) + O(m^{-2}).$$
Thus
$$\mu_{mc}(\I_{\PP(E|_{B'})},L_m) - \mu(\PP(E),L_m) = \frac{1}{m}[\mu(\I_{B'},\O_B(1))-\mu(B)] + O(m^{-2}),$$
as required.
\end{proof}\vskip 5pt

\begin{proof}[Proof of Theorem \ref{thm:slopestableprojectivebundles}]  
  By Proposition \ref{prop:seshadri} there is an $m_0$ such that for
  all $m\ge m_0$, $\epsilon(\PP(F),L_m)=1$ for all saturated coherent
  subsheaves $F<E$ with $\mu\_F\ge \mu\_E$.
  
  As the family of destabilising subsheaves of $F$ of $E$ is
  bound\-ed, the set $\{c_i(F)\in H^{2i}(B)\colon F<E,
  \mu\_F\ge\mu\_E, 0\le i\le n \}$ is finite. Thus we can bound the
  $O(\widetilde{m}^{2b-3})$ terms in Proposition \ref{prop:slopeprojectivebundles}
  independently of $F$.  Furthermore there is a $\delta>0$ (again
  independent of $F$) such that $\mu\_F> \mu\_E$ implies $\mu\_F\ge
  \mu\_E+\delta$.  Hence for all saturated coherent subsheaves $F<E$
  with $\mu\_F\ge \mu\_E$ and $m\ge m_0$,
\begin{eqnarray} 
 &&\hspace{-6mm}\mu_1(\O_{\PP(F)},L_m)-\mu(\PP(E),L_m)\nonumber \\
&=& C\left((r+1)(\mu\_E-\mu\_F)\widetilde{m}^{2b-1} + O(\widetilde{m}^{2b-2})\right)\nonumber\\
&\le& -C(r+1)\big(\delta m^{2b-1} + O(m^{2b-2})\big), \label{Cdelta}
\end{eqnarray}
where $C=C(\widetilde{m})>0$ is independent of $F$.

Now suppose that $E$ is not slope semistable.  Then there exists a
coherent $F<E$ with $\mu\_F>\mu\_E$.  Replace $F$ by its saturation
(\emph{i.e.}\ the kernel of $E\to (E/F)/$torsion), which has slope
$\ge\mu\_F>\mu\_E$.  Making $m_0$ larger if necessary we have that
$\mu_1(\O_{\PP(F)},L_m)<\mu(\PP(E),L_m)$ for $m\ge m_0$ by
(\ref{Cdelta}), so $(\PP(E),L_m)$ is not slope semistable.
 
  Similarly if $(B,\O_B(1))$ is not slope semistable then there is a
  $B'$ and $c$ with $c<\epsilon(B)$ and $\mu_c(\I_B)>\mu(B)$.
  Therefore, for $m\gg 0$, $c+O(m^0)/m<\epsilon(B)$, so by Proposition \ref{b'},
  $mc<\epsilon(\PP(E|_{B'}),L_m)$ and $\mu_{cm}(\I_{\PP(E|_{B'})},L_m)
  > \mu(\PP(E),L_m)$.  Thus $\PP(E|_{B'})$ strictly destabilises
  $(\PP(E),L_m)$ for $m\gg 0$.
\end{proof}\vskip 5pt

\subsection{Unstable blow ups}\label{sec:blowups}\ \vskip 5pt
\noindent
Fix $Z\subset X$. If we form the blow up $\pi\colon\widehat X\to X$ of
$X$ along $Z$, with exceptional divisor $E$, then since for $k\gg0$
$$
H^0_X(L^{\otimes k}\otimes\I_{\!Z}^{xk})\cong H^0_{\widehat
  X}(\pi^*L^{\otimes k}\otimes\I_E^{xk}),
$$
there is a strong link between $Z\subset X$ and $E\subset\widehat
X$. Morally, $Z$ destabilises $(X,L)$ if and only if $E$ destabilises
$(\widehat X,\pi^*L)$, but the latter line bundle is only semi-ample.
However, $L_d:=\pi^*L(-dE)$ is ample for $0<d<\epsilon(Z)$, and with
respect to this polarisation the Seshadri constants are related by
$\epsilon(E)=\epsilon(Z)-d$. For $k\gg0$,
$$
H^0_{\widehat X}(L_d^{\otimes k}\otimes\I_E^{xk})\cong
H^0_X(L^{\otimes k}\otimes\I_{\!Z}^{(x+d)k}),
$$
so for $d<c\le \epsilon(Z)$,
$$
\mu_{c-d}(\I_E,L_d) = \frac{\int_d^c \left(a_1(x) +
  \frac{a_0'(x)}{2}\right)dx}{\int_d^c a_0(x)dx}\,.
$$
As $d\to0$, $\mu_{c-d}(\I_E)\to\mu_c(\I_Z)$ and $\mu(\widehat
X)\to\mu(X)$ as expected.

This can be applied in the following way. Suppose that a singular
point strictly destabilises a variety $X$, and that its blow up is
smooth. More generally, fix an ideal sheaf $\I_Z\subset\O_X$ whose
blow up is smooth (this exists by resolution of singularities) and
suppose that $\I_Z$ strictly destabilises $(X,L)$. Then for small $d$,
$(\widehat X,L_d)$ is also strictly unstable, and so has no cscK
metric.  This gives an easy way of producing smooth polarised
varieties without cscK metrics.  

\subsection{Unstable rational manifolds}\label{sec:rational} \vskip 5pt
\begin{example}\label{p2blownupat1point}
  $\PP^2$ \textbf{blown up at 1 point}. Any polarisation on $\pi\colon
  X\to\PP^2$ is a multiple of $L=L_q = \O_{\PP^2}(1) - qE$ for
  $q\in(0,1)\cap\Q$.  We claim that the exceptional curve $E$
  destabilises $(X,L)$ for any such $q$.  Let $Z=E$.  Then
  $L(-cZ)=\O_{\PP^2}(1)\big(-(q+c)E\big)$ is nef for $c<1-q$, hence
  $\epsilon=\epsilon(Z,L)=1-q$.  By (\ref{curvsurf})
\begin{eqnarray*}
\mu(X,L) &=& \frac{3-q}{1-q^2}\,, \quad \text{and}\\
\mu_\epsilon(\O_Z)  &=& \frac{3}{(1-q)(2q+1)}\,.
\end{eqnarray*}
Since
$$(3-q)(1-q)(2q+1)-3(1-q^2)=2q(1-q)^2>0,$$
we have
$\mu_\epsilon(\O_Z)<\mu(X,L)$ for all $0<q<1$.

In fact this example is covered by Section \ref{proj}, since $X$ is a
projective bundle $X\cong\PP(\O_{\PP^1}\oplus\O_{\PP^1}(1))$ (with $E$
the projectivisation of the destabilising subbundle
$\O_{\PP^1}(1)$). It is significant that $E$ destabilises
$X$ for \emph{all} $q$, since \emph{a priori} $X$ is unstable
for all polarisations because of the non-reductivity of Aut($X$).
\end{example}

To find destabilising examples it is convenient to allow $L$ to tend
to a divisor that is not necessarily ample.

By analogy with
(\ref{definitionofa0x}, \ref{definitionofa1x}) we use the Riemann-Roch
formula on $\widehat X$ to define the slope with respect to any
divisor $F$ by
\begin{eqnarray*}
  a_0^F(x) &=& \frac{1}{n!}\int_{\widehat X} c_1(F(-xE))^n,\nonumber\\
a_1^F(x) &=& -\frac{1}{2(n-1)!}\int_{\widehat X} c_1(K_{\widehat X}).c_1(F(-xE))^{n-1},
\end{eqnarray*}
and (\ref{qslope})
\begin{eqnarray*}
 \mu(X,F) &=& -\frac{n\int_X c_1(K_X).c_1(F)^{n-1}}{2\int_X c_1(F)^n},\\
 \mu_c(\O_Z,F)&=&\frac{\int_0^c \left(\tilde{a}^F_1(x) + \frac{\tilde{a}^F_0\,\!' (x)}{2}\right)
dx}{\int_0^c \tilde{a}^F_0(x)dx}\,,
\end{eqnarray*}
where $\tilde{a}^F_i(x) = a^F_i(0)-a^F_i(x)$.  Note that since $F$ is not
assumed to be ample, these could be infinite.
\begin{prop}\label{destabilisingusingnefdivisors}
  Let $F$ be a nef divisor on $X$.  Suppose that there is a $c>0$
  such that $F(-cE)$ is nef on $\widehat X$, and
  $$
  \int_X c_1(F)^n  \int_0^c \left(\tilde{a}^F_1(x) + \frac{\tilde{a}_0^F\,\!'(x)}{2}\right)
  dx$$
is strictly less than
$$ - \frac{n\int_X c_1(K_X).c_1(F)^{n-1}}{2} \int_0^c \tilde{a}^F_0(x) dx.$$
  (In particular this holds if $\int_X c_1(F)^n>0$ and
  $$\mu_c(\O_Z,F)<\mu(X,F)<\infty).$$

  Then $Z$ strictly destabilises $(X,L)$ for $L$ sufficiently close to
  $F$.  More precisely: if $G$ is an ample divisor and $L=F(\delta G)$, then there
  is a $\delta_0>0$ such that $Z$ strictly destabilises $(X,L)$ for
  all $0<\delta<\delta_0$.
\end{prop}
\begin{proof}
  Since $L(-cE)=F(-cE+\delta G)$ is nef we have $\epsilon(Z,L)\ge c$.
  Notice that $ \int_X c_1(L)^n \int_0^c \tilde{a}_0^L(x) dx
  [\mu(\O_Z,L)-\mu(X,L)]$ equals 
  $$\int_X \!\! c_1(L)^n \int_0^c \!\left(\tilde{a}^L_1(x) + \frac{\tilde{a}_0^L\,\!'(x)}{2}
 \right) dx + \frac{n\int_X c_1(K_X).c_1(L)^{n-1}}{2} \int_0^c \!\tilde{a}^L_0(x) dx.$$
  As
  $\delta$ tends to zero this tends to
  $$
  \int_X c_1(F)^n \!\! \int_0^c \!\left(\tilde{a}^F_1(x) + \frac{\tilde{a}_0^F\,\!'(x)}{2}\right)
  dx + \frac{n\int_X c_1(K_X).c_1(F)^{n-1}}{2} \int_0^c \!\tilde{a}^F_0(x) dx,$$
  which
  is assumed to be strictly negative.  Since $L$ is ample, and
  $c\le\epsilon(Z,L)$, $\int_X c_1(L)^n\int_0^c \tilde{a}_0^L(x)dx>0$.  Thus
  $\mu_c(\O_Z,L)<\mu(X,L)$ for $\delta$ sufficiently small.
\end{proof}

\begin{cor}{\bf Unstable blow ups}. Suppose
that $(X,L)$ is destabilised by $Z$, and let $Y$ be the blow up of $X$ along
a centre disjoint from $Z$.  Then for polarisations making the exceptional
set small, $Y$ is destabilised by the proper transform of $Z$.
\end{cor}

\begin{example}
  $\PP^2$ \textbf{blown up at m distinct points}.  Let $X$ be $\PP^2$
  blown up at $m\ge 1$ distinct points, with exceptional divisors
  $\{E_i\}_{i=1}^m$. Then applying the above to Example
  \ref{p2blownupat1point} shows that $X$ is slope
  unstable with respect to suitable polarisations: those of the form 
  $\O_{\PP^2}(1)\big(-\sum_{i=1}^mq_iE_i\big)$ with $0<q_i\ll q_1<1$ for $i\ge2$.
  \end{example}

\begin{rmk}
It is important to note that these polarisations are far from the anticanonical
polarisation. For
generic configurations of points $K_X^*$ is ample if
$m\le 8$, and this polarisation \emph{does} admit a cscK (in fact
K\"ahler-Einstein) metric, unless $m=1$ or $m=2$ \cite{Ti1}.
\end{rmk}  

\begin{rmk}\label{rmk:folklore}{\bf The folklore conjecture}.
The case of $\PP^2$ blown up at $\ge4$ points gives smooth polarised
  del Pezzo surfaces with $\aut(X)=0$ but no cscK metric in certain
  classes. This is in contrast to the case of the anticanonical
  polarisation for which Tian \cite{Ti1} proved the ``folklore
  conjecture", that smooth Fano surfaces have a K\"ahler-Einstein metric
  if and only
  if their holomorphic automorphism group is reductive, and disproves
  the conjecture for manifolds of dimension $n\ge 3$.  There are also
  examples of ruled surfaces \cite{BdB} which show the folklore
  conjecture for cscK metrics on surfaces does not hold.
\end{rmk}

\begin{rmk}{\bf Unstable elliptic surface}. If $X$ is $\PP^2$ blown up
  at 9 points which are the intersection of two cubics, then $X$ is
  slope unstable with respect to suitable polarisations.  Thus we have
  a polarised elliptic surface (``half a $K3$ surface") which is not
  K-semistable and does not admit a cscK metric.
\end{rmk}

\begin{example}\textbf{$-2$ curves}. We now give an example of a destabilising
  $-2$-curve on a del Pezzo surface. 
  Blow $\PP^2$ up at a point, and let $X$ be its blow up at a
  point on the exceptional divisor.  Thus $X$ contains a $-2$ curve
  $E_1$ and an exceptional $-1$-curve $E_2$.  Thus $E_1^2=-2$,
  $E_2^2=-1$ and $E_1.E_2=1$.  Notice that
  $\O_{\PP^2}(1)\big(\!-\frac12E_1-E_2\big)$ is nef ($X$ is toric,
  and the line bundle's degree on each invariant curve is nonnegative, so
  the toric Kleiman criterion applies).  Also $\O_{\PP^2}(1)$ is nef, as
  is $\O_{\PP^2}(1)\big(-E_1-E_2\big)$ since it is the pullback of a nef bundle on
  $\PP^2$ blown up at one point.
  By convexity of the ample cone, $\O_{\PP^2}(1)\big(\!-qE_1 -r E_2\big)$ is ample
  for $1\ge r\ge q\ge r/2\ge0$.
 
  Set $L=\O_{\PP^2}(1)\big(\!-\frac12E_1-rE_2\big)$, which is ample for $1>r>1/2$,
  and let $Z=E_1$ be the exceptional $-2$-curve in $X$.  From the
  above, $c_r=r-\frac12\le \epsilon(Z,L)$.  Then $L.Z=1-r \rightarrow
  0$ as $r\rightarrow 1$, so (\ref{curvsurf})
\begin{eqnarray*}
 \mu(X,L) &=& \frac{\O_{\PP^2}(3)\big(\!-E_1-2E_2\big).L}{L^2}=\frac{6-2r}{1+2r-2r^2} \rightarrow 4  \,\, \text{as $r\rightarrow 1$}, \\
 \mu_{c_r}(\O_Z,L) &=& \frac{3(L.Z+c_r)}{c_r(L.Z+2c_r)}\rightarrow 3 \,\, \text{as $r\rightarrow 1$}.
\end{eqnarray*}
Thus for $r$ close to $1$, $\mu_{c_r}(\O_Z)<\mu(X,L)$ so $Z$ strictly
destabilises $(X,L)$. Further blowing up $X$ in some points disjoint from
the $-2$-curve and taking a polarisation in which the new exceptional divisors
are small, this also gives examples with no automorphism group.
\end{example} 

\begin{example}{\bf $\PP^n$ blown up at points}. The exceptional divisor
$E$ strictly destabilises $\PP^n$ blown up at a point with respect to all
polarisations. This is an application of
Theorem \ref{divisors}; we omit the gory details.

Therefore the blow up of $\PP^n$ at $m\ge 1$ distinct points (or a point and
some disjoint subvarieties) is
  unstable with respect to polarisations which make one component of
  the exceptional set large, and the other $m-1$ small.
\end{example}

\vskip 4mm

{\small \noindent {\tt jaross@math.columbia.edu} \\
\noindent \small{\tt richard.thomas@imperial.ac.uk}} \newline
\noindent Department of Mathematics, Columbia University, \\New York, NY 10027.
USA. \\
Department of Mathematics, Imperial College, \\
London SW7 2AZ. UK.

\end{document}

%% file: cone.pstex_t
\begin{picture}(0,0)%
\includegraphics{cone.pstex}%
\end{picture}%
\setlength{\unitlength}{4144sp}%
\begingroup\makeatletter\ifx\SetFigFont\undefined%
\gdef\SetFigFont#1#2#3#4#5{%
  \reset@font\fontsize{#1}{#2pt}%
  \fontfamily{#3}\fontseries{#4}\fontshape{#5}%
  \selectfont}%
\fi\endgroup%
\begin{picture}(3577,3254)(1883,-3575)
\put(2693,-916){\makebox(0,0)[lb]{\smash{\SetFigFont{12}{14.4}{\rmdefault}{\mddefault}{\updefault}{\color[rgb]{0,0,0}$X_t$}%
}}}
\put(3829,-487){\makebox(0,0)[lb]{\smash{\SetFigFont{12}{14.4}{\rmdefault}{\mddefault}{\updefault}{\color[rgb]{0,0,0}$Z'$}%
}}}
\put(4037,-1870){\makebox(0,0)[lb]{\smash{\SetFigFont{9}{10.8}{\rmdefault}{\mddefault}{\updefault}{\color[rgb]{0,0,0}$E$}%
}}}
\put(1883,-1605){\makebox(0,0)[lb]{\smash{\SetFigFont{12}{14.4}{\rmdefault}{\mddefault}{\updefault}{\color[rgb]{0,0,0}$Z\times\{t\}$}%
}}}
\put(3380,-3526){\makebox(0,0)[lb]{\smash{\SetFigFont{12}{14.4}{\rmdefault}{\mddefault}{\updefault}{\color[rgb]{0,0,0}$\X_0=\widehat X\cup_EP$}%
}}}
\end{picture}

%% file: toriccone.pstex_t
\begin{picture}(0,0)%
\includegraphics{toriccone.pstex}%
\end{picture}%
\setlength{\unitlength}{1579sp}%
\begingroup\makeatletter\ifx\SetFigFont\undefined%
\gdef\SetFigFont#1#2#3#4#5{%
  \reset@font\fontsize{#1}{#2pt}%
  \fontfamily{#3}\fontseries{#4}\fontshape{#5}%
  \selectfont}%
\fi\endgroup%
\begin{picture}(6675,7862)(1426,-8061)
\put(3001,-6961){\makebox(0,0)[lb]{\smash{\SetFigFont{7}{8.4}{\rmdefault}{\mddefault}{\updefault}{\color[rgb]{0,0,0}$f=c$}%
}}}
\put(5476,-7036){\makebox(0,0)[lb]{\smash{\SetFigFont{7}{8.4}{\rmdefault}{\mddefault}{\updefault}{\color[rgb]{0,0,0}$f=\sum_{i=1}^m\frac{f_i}{m_i}$}%
}}}
\put(1426,-1261){\makebox(0,0)[lb]{\smash{\SetFigFont{8}{9.6}{\rmdefault}{\mddefault}{\updefault}{\color[rgb]{0,0,0}$Q$}%
}}}
\put(3076,-1561){\makebox(0,0)[lb]{\smash{\SetFigFont{7}{8.4}{\rmdefault}{\mddefault}{\updefault}{\color[rgb]{0,0,0}$P_c$}%
}}}
\put(3526,-4036){\makebox(0,0)[lb]{\smash{\SetFigFont{7}{8.4}{\rmdefault}{\mddefault}{\updefault}{\color[rgb]{0,0,0}$P_0$}%
}}}
\put(3601,-5611){\makebox(0,0)[lb]{\smash{\SetFigFont{7}{8.4}{\rmdefault}{\mddefault}{\updefault}{\color[rgb]{0,0,0}$p$}%
}}}
\put(8101,-4786){\makebox(0,0)[lb]{\smash{\SetFigFont{7}{8.4}{\rmdefault}{\mddefault}{\updefault}{\color[rgb]{0,0,0}$0$}%
}}}
\put(6076,-7936){\makebox(0,0)[lb]{\smash{\SetFigFont{7}{8.4}{\rmdefault}{\mddefault}{\updefault}{\color[rgb]{0,0,0}$Z$}%
}}}
\put(6451,-2461){\makebox(0,0)[lb]{\smash{\SetFigFont{7}{8.4}{\rmdefault}{\mddefault}{\updefault}{\color[rgb]{0,0,0}$\pi$}%
}}}
\put(1501,-6211){\makebox(0,0)[lb]{\smash{\SetFigFont{7}{8.4}{\rmdefault}{\mddefault}{\updefault}{\color[rgb]{0,0,0}$P$}%
}}}
\put(8101,-2236){\makebox(0,0)[lb]{\smash{\SetFigFont{7}{8.4}{\rmdefault}{\mddefault}{\updefault}{\color[rgb]{0,0,0}$c$}%
}}}
\end{picture}

%% file: slopeobs_archive.bbl
\begin{thebibliography}{ATGT-F}
\bibitem[AT]{ATF} Apostolov, V. and T{\o}nnesen-Friedman, C. (2004).
  \emph{A remark on K\"ahler metrics of constant scalar curvature on
    ruled complex surfaces.}  Preprint math.DG/0411271.
  
\bibitem[ACGT]{ACGT} Apostolov, V., Calderbank, D., Gauduchon, P. and
  T{\o}nnesen-Friedman, C. (2005).  \emph{Hamiltonian 2-forms in Kahler
    geometry, III Extremal metrics and stability.} Preprint
  math.DG/0511118.


\bibitem[Au]{Au} Aubin, T. (1976).  \emph{\'Equations du type
    Monge-Amp\`ere sur les vari\'et\'es k\"ahleriennes compactes.} C.
  R. Acad. Sci. Paris S\'er. A-B \textbf{283}, A119--A121, MR0433520, Zbl 0333.53040.
  
\bibitem[BdB]{BdB} Burns, D. and De Bartolomeis, P. (1988).
  \emph{Stability of vector bundles and extremal metrics.} Invent.
  Math. \textbf{92}, 403--407, MR0936089, Zbl 0645.53037.
  
\bibitem[LeB]{LeB} LeBrun, C. (1995).  \emph{Polarized $4$-manifolds,
    extremal K\"ahler metrics, and Seiberg-Witten theory.}  Math. Res.
  Lett. \textbf{2}, 653--662, MR1359969, Zbl 0874.53051.

\bibitem[LeBS]{LeBS} LeBrun, C. and Simanca, S. (1994).
  \emph{Extremal K\"ahler metrics and complex deformation theory.}
  GAFA \textbf{4}, 298--336, MR1274118, Zbl 0801.5305.

\bibitem[CEL]{CEL} Cutkosky, S., Ein, L. and Lazarsfeld, R. (2001).
  \emph{Positivity and complexity of ideal sheaves.}  Math. Ann.
  \textbf{321}, 213--234, MR1866486, Zbl 1029.14022.
  
\bibitem[CT]{CT} Chen, X. and Tian, G. (2003).  \emph{Geometry of
    K\"ahler metrics and holomorphic foliation by discs.}  Preprint
  math.DG/0409433.
  
\bibitem[DP]{DP} Demailly, J.-P. and Paun, M. (2004).  \emph{Numerical
    characterization of the K\"ahler cone of a compact K\"ahler
    manifold}.  Ann. of Math. (2)  \textbf{159}, 1247--1274, MR2113021, Zbl 1064.32019.

  
\bibitem[Do1]{Do1} Donaldson, S.  (1997).  \emph{Remarks on gauge
    theory, complex geometry and $4$-manifold topology}.  Fields
  Medallists' lectures, 384--403, World Sci. Publishing, MR1622931.
  
\bibitem[Do2]{Do2} Donaldson, S. (2001).  \emph{Scalar curvature
    and projective embeddings, I.}  Jour. Diff. Geom. \textbf{59},
  479--522, MR1916953, Zbl 1052.32017.
  
\bibitem[Do3]{Do3} Donaldson, S. (2002).  \emph{Scalar curvature
    and stability of toric varieties}. Jour. Diff. Geom.  \textbf{62},
  289--349, MR1988506, Zbl pre02171919.
  
\bibitem[Do4]{Do4} Donaldson, S. (2005). \emph{Scalar curvature and
    projective embeddings II}.
  Quarterly Journal of Mathematics \textbf{56}, 345--356, MR2161248.
  
\bibitem[Do5]{Do5} Donaldson, S. (2005). \emph{Lower bounds on the
    Calabi Functional}. Preprint math.DG/0506501.

\bibitem[Ful]{Ful}
Fulton, W. (1977).
\emph{A Hirzebruch-Riemann-Roch formula for analytic spaces and non-projective
algebraic varieties.} Compositio Math. \textbf{34}, 279--283, MR0460323, Zbl 0367.14008.

  
\bibitem[Fut]{Fut} Futaki, A. (1983).  \emph{On compact K\"ahler
  manifolds of constant scalar curvatures.}  Proc. Japan Acad. Ser. A
Math. Sci. \textbf{59}, 401--402, MR0726535, Zbl 0539.53048.

\bibitem[HL]{HL}
Huybrechts, D. and Lehn, M. (1997).
\emph{Geometry of moduli spaces of shaves.}
Aspects in Mathematics Vol. E31, Vieweg, MR1450870, Zbl 0872.14002.

\bibitem[Ha]{Ha} Hartshorne, R. (1977).  \emph{Algebraic Geometry.}
  Graduate Texts in Mathematics 52, Springer-Verlag, MR0463157, Zbl 0531.14001.
  
\bibitem[Ho]{Ho} Hong, Y-J. (1999).  \emph{Constant Hermitian scalar
    curvature equations on ruled manifolds}, Jour. Diff. Geom.
  \textbf{53}, 465--516, MR1806068, Zbl pre01782630.

\bibitem[Kl]{Kl} Kleiman, S. L. (1966).  \emph{Toward a numerical
    theory of ampleness.} Ann. of Math. \textbf{84}, 293--344, MR0206009, Zbl 0146.17001.

\bibitem[La]{La}
Lazarsfeld, R. (2003).
\emph{Positivity in algebraic geometry. II. Positivity for vector bundles, and multiplier ideals.}
Ergeb. Math. Grenzgeb. (3), Springer-Verlag, MR2095472, Zbl pre02134815.
 
\bibitem[Mab]{Mab} Mabuchi, T. (2004).  \emph{An energy-theoretic
    approach to the Hitchin-Kobayashi correspondence for manifolds,
    I.} Invent. Math. \textbf{159},  225--243, MR2116275, Zbl pre02156015.
  
\bibitem[Mo]{Mo} Morrison, I. (1980).  \emph{Projective stability of
    ruled surfaces.}  Invent. Math. \textbf{56}, 269--304, MR0561975, Zbl 0423.14005.
  
\bibitem[Mu]{Mu} Mumford, D. (1977).  \emph{Stability of projective
    varieties.}  Enseignement Math. (2) \textbf{23}, 39--110, MR0450273, Zbl 0376.14007.
  
\bibitem[GIT]{GIT} Mumford, D., Fogarty, J. and Kirwan, F. (1994).
  \emph{Geometric Invariant Theory.}  Third edition, Erg. Math.
  \textbf{34}, Springer-Verlag, Berlin, MR1304906, Zbl 0797.14004.

\bibitem[Na]{Na}
Nadel, A. (1990).
\emph{Multiplier ideal sheaves and K\"ahler-Einstein metrics of positive scalar
curvature.} Ann. of Math. \textbf{132}, 549--596, MR1078269, Zbl 0731.53063.
  
\bibitem[P]{P} Paul, S. (2004).  \emph{Geometric analysis of Chow
    Mumford stability.} Adv. Math. \textbf{182}, 333-356, MR2032032, Zbl 1050.53061.
  
\bibitem[PT]{PT} Paul, S. and Tian, G. (2004).  \emph{Algebraic and
    Analytic K-Stability.} Preprint math.DG/0405530.

\bibitem[Ro]{Ro}
Ross, J. (2003).
\emph{Instability of polarised algebraic varieties.} PhD thesis, Imperial
College.
  
\bibitem[RT]{RT} Ross, J. and Thomas, R. (2004) \emph{A study of the Hilbert-Mumford
criterion for the stability of projective varieties.} Preprint math.AG/0412519.

\bibitem[Ti1]{Ti1} Tian, G. (1990).  On Calabi's conjecture for
  complex surfaces with positive first Chern class.  Invent. Math.
  \textbf{101}, 101--172, MR1055713, Zbl 0716.32019.
  
\bibitem[Ti2]{Ti2} Tian, G. (1994).  The $K$-energy on hypersurfaces
  and stability.  Comm. Anal. Geom. \textbf{2}, 239--265, MR1312688, Zbl 0846.32019.
  
\bibitem[Ti3]{Ti3} Tian, G. (1997).  \emph{K\"ahler-Einstein metrics
    with positive scalar curvature.}  Invent. Math. \textbf{130},
  1--37, MR1471884, Zbl 0892.53027.
  
\bibitem[Wa]{Wa} Wang, X. (2004) \emph{Moment map, Futaki invariant and
    stability of projective manifolds.} Comm. Anal. Geom. \textbf{12}, 1009--1037, MR2103309, Zbl pre02147072.


\bibitem[We]{We} Weinkove, B. (2004).  \emph{On the J-flow in higher
    dimensions and the lower boundedness of the Mabuchi energy.}
  Preprint math.DG/0309404.

\bibitem[Y1]{Y1} Yau, S.-T. (1978).  \emph{On the Ricci curvature of a
    compact K\"ahler manifold and the complex Monge-Amp\`ere equation.
    I.} Comm. Pure Appl. Math. \textbf{31}, 339--411, MR0480350, Zbl 0369.53059.
  
\bibitem[Y2]{Y2} Yau, S.-T. (1993).  \emph{Open problems in geometry.}
  Differential geometry: partial differential equations on manifolds
  (Los Angeles, CA, 1990), 1--28, Proc. Sympos. Pure Math.
  \textbf{54}, AMS publications, MR1216573.
  
\bibitem[Zh]{Zh} Zhang, S. (1996).  \emph{Heights and reductions of
    semi-stable varieties.}  Compositio Math. \textbf{104}, 77--105,
  MR1420712, Zbl 0924.11055.

\end{thebibliography}
